\documentclass[10pt,twocolumn,letterpaper]{article}

\usepackage{graphicx}
\usepackage{amsmath}
\usepackage{amssymb}
\usepackage{booktabs}

\usepackage[pagebackref,breaklinks,colorlinks]{hyperref}

\usepackage{enumitem}

\newif\ifTWOCOLUMNS
\TWOCOLUMNStrue

\DeclareMathOperator*{\argmin}{arg\,min}
\DeclareMathOperator*{\argmax}{arg\,max}
\DeclareMathOperator*{\dom}{dom}

\def\myparagraph#1{\vspace{2pt}\noindent{\bf #1~~}}

\newcommand{\eqdef}{{\stackrel{\mbox{\tiny \tt ~def~}}{=}}}

\def\calW{{\cal W}}

\long\def\ignore#1{}
\def\myps[#1]#2{\includegraphics[#1]{#2}}

\def\br(#1,#2){{\langle #1,#2 \rangle}}
\def\setZ[#1,#2]{{[ #1 .. #2 ]}}

\def\true{\mbox{\tt true}}
\def\false{\mbox{\tt false}}

\def\calX{\dot{\cal X}}

\def\X{{\cal X}}
\def\Y{{\cal Y}}

\def\conv{\operatorname{conv}}

\def\q={\quad=\quad}
\def\qq={\qquad=\qquad}

\def\calE{{\cal E}}

\def\calG{{\cal G}}

\def\calL{{\cal L}}

\def\calS{{\cal S}}

\def\calV{{\cal V}}

\def\psfile[#1]#2{}
\def\psfilehere[#1]#2{}

\def\assign(#1,#2){\langle#1,#2\rangle}
\def\edge(#1,#2){(#1,#2)}

\def\slack(#1){\texttt{slack}({#1})}
\def\barslack(#1){\overline{\texttt{slack}}({#1})}

\def\unitvec(#1){{{\bf u}_{#1}}}

\date{}

\begin{document}

\title{Solving relaxations of MAP-MRF problems: \\ Combinatorial in-face Frank-Wolfe directions}

\author{Vladimir Kolmogorov\\
Institute of Science and Technology Austria (ISTA) \\
Am Campus 1, Klosterneuburg 3400, Austria \\
{\tt\small vnk@ist.ac.at\vspace{-8pt}}
}

\maketitle

\begin{abstract}
  We consider the problem of solving LP relaxations of MAP-MRF inference problems,
  and in particular the method proposed recently in~\cite{Swoboda:CVPR19,KP:ICML21}.
  As a key computational subroutine, it uses a variant of the Frank-Wolfe (FW) method
  to minimize a smooth convex function over a combinatorial polytope.
  We propose an efficient implementation of this subroutine based
  on {\em in-face Frank-Wolfe directions}, introduced in~\cite{FW:in-face} in a different context.
  More generally, we define an abstract data structure for a combinatorial subproblem
   that enables in-face FW directions, and describe its specialization
  for tree-structured MAP-MRF inference subproblems.
  Experimental results indicate that the resulting method is the current state-of-art LP solver
  for some classes of problems. Our code is available at \url{pub.ist.ac.at/~vnk/papers/IN-FACE-FW.html}.
\end{abstract}

\section{Introduction}\label{sec:intro}
The main focus of this paper is on the problem of minimizing a function of discrete variables $z=(z_1,\ldots,z_n)$ with unary and pairwise terms:
	\begin{equation}\label{eq:MAPMRF}
	  \min_{z\in D_1\times\ldots \times D_n} \;\; \sum_{v\in [n]} f_v(z_v)+ \sum_{uv\in\calE} f_{uv}(z_u,z_v)
	\end{equation}
Here $\calG=([n],\calE)$ is an undirected graph and $D_1,\ldots,D_n$ are finite sets.
This problem is often referred to as {\em MAP-MRF inference} ({\em maximum a posteriori inference in a Markov Random Field}).

A prominent approach to tackle this NP-hard problem in practice is to solve its natural LP relaxation (see e.g.~\cite{werner2007linear}),
also called {\em Basic LP relaxation}~\cite{kolmogorov15:power}:
\begin{subequations}\label{eq:BLP}
\begin{align}
\hspace{0pt} &  \min_{\xi\ge {\bf 0}} \;\; \sum_{\substack{v\in [n] \\ a\in D_v}} f_v(a)\xi_{va} + \!\!\!\!\!\!\sum_{\substack{uv\in\calE \\ (a,b)\in D_u\times D_v}} \!\!\!\!\!\! f_{uv}(a,b)\xi_{ua;vb} \\
\hspace{0pt} & \sum_{b'\in D_v}\! \xi_{ua;vb'} = \xi_{ua}, 
 \sum_{a'\in D_u}\!\xi_{ua';vb} = \xi_{vb}  \quad \forall uv,a,b  \\
\hspace{0pt} & \sum_{a\in D_v} \xi_{ua}=1 \quad\forall v
\end{align}
\end{subequations}
Designing algorithms to (approximately) solve this relaxation for large-scale problems has been a very active area of research.
A popular approach is to use {\em message passing} techniques,
which perform a block-coordinate ascent on the dual objective~\cite{Kolmogorov2006,Werner07,MPLP,SRMPKolmogorov,IRPSLP,Tourani:AISTATS20}
This strategy is very effective for some problems, but for other problems it may get stuck in a suboptimal point.
Many techniques have been developed that are guaranteed to converge to the optimal solution of the LP relaxation
\cite{
StorvikDahlLagrangeanBasedMAP, 
SchlesingerSubgradient, 
LagrangeanRelaxationJohnsonMalioutov, 
RavikumarProximalMethodsMAPMRF, 
AcceleratedMAPJojic, 
savchynskyy2011study, 
schmidt2011evaluationProximalMAP, 
DualDecompositionKomodakis, 
Martins:ICML11,
AdaptiveDiminishingSmoothingSavchynskyyUAI, 
MapMirrorDescent, 
Schwing:NIPS12,
Schwing:ICML14,
Swoboda:CVPR19,
KP:ICML21}. 

In this paper we revisit the approach in~\cite{Swoboda:CVPR19,KP:ICML21}. Its key computational subroutine is
to minimize a quadratic convex function over combinatorial polytope, which is done by invoking 
a variant of the Frank-Wolfe (FW) algorithm~\cite{FWolfe56}. We study efficient implementations of the latter
in the context of MAP-MRF inference. Our main contribution is incorporating
{\em in-face FW directions}
introduced in~\cite{FW:in-face}. The  idea is to speed-up computations by running FW algorithm on a smaller ``contracted''
subproblem obtained by taking a face of the polytope containing the current point. 
It has been used for applications such as low-rank matrix completion~\cite{FW:in-face},
cluster detection in networks~\cite{sdefectiveclique}, and training sparse neural networks with $\ell_1$ regularization~\cite{Grigas:19}.
We investigate the use of in-face FW directions for general combinatorial polytopes, and describe
an abstract data structure that enables such directions. We then specialize it to subproblems corresponding to tree-structured MAP-MRF inference problems.
Our approach has the following features:

\begin{itemize}[topsep=0pt, itemsep=-3pt, leftmargin=11pt]
\item It may happen that the contracted subproblem splits into independent subproblems. These subproblems are handled by a block-coordinate version of FW.  
\item We store a cache of ``atoms'' for each contracted subproblem. We describe how to efficiently transform these atoms when the current face is recomputed. 
\item For an edge $uv\in \calE$ and fixed fractional unary vectors for $u,v$
we can compute an optimal fractional pairwise vector for edge $uv$ by
solving a small-scale optimal transportation (OT) problems. Such computations were used in~\cite{SavchynskyySchmidt:13} for computing primal feasible solutions of  relaxation~\eqref{eq:BLP}.
We show how to use them for improving the performance of in-face FW directions.
\end{itemize}
To our knowledge, the issues above have not been discussed in the literature so far.

We remark that in-face FW directions effectively implement the following rather natural idea:
the optimization should be performed only over ``active'' pairs $(v,a)$ that are likely to
be present in the support of an optimal solution; the other pairs should be fixed.
A related idea appeared in the context of message passing algorithms in~\cite{Tarlow:ICML11},
where messages are updated only in a subgraph in which the current best labels keep changing.
The method in~\cite{Tarlow:ICML11} works only with dual variables,
and uses heuristic criteria for choosing the subgraph.
We believe that in-face FW directions allow a more principled criterion for deciding which variables to fix and for how long.
Note that Freund et al.~\cite{FW:in-face} proved that their criterion retains the convergence rate of the basic FW algorithm.
We use a different criterion that also takes into account the ratio of runtimes on the original and on contracted subproblems.

In Section~\ref{sec:experiments} we test the algorithms on benchmark problems in the evaluation~\cite{OpenGMBenchmark},
and compare them with LP solvers used in~\cite{OpenGMBenchmark}.
Results suggest that the method in~\cite{Swoboda:CVPR19,KP:ICML21} with in-face FW directions is the current state-of-the-art LP solver for certain classes
of problems. 

\section{Background}
We will consider a more general problem of minimizing a function represented as a sum of {\em tractable subproblems}:
\begin{equation}
  \label{eq:decomposition}
  \min_{X \in \mathbb R^d} f(X), \quad f(X) := \sum_{t\in T} f_t(X_{A_t})
\end{equation}
Here term $t\in T$ is specified by a subset of variables $A_t\subseteq[d]$
and a function $f_t:\mathbb R^{A_t}\rightarrow\mathbb R\cup\{+\infty\}$ of $|A_t|$ variables.
Vector $X_{A_t}\in\mathbb R^{A_t}$ is the restriction of vector $X\in\mathbb R^d$ to $A_t$.
We assume that the effective domain $\dom f_t=\{X\in\mathbb R^{A_t}\:|\:f_t(X)<+\infty\}$
is a finite non-empty set (implying that~\eqref{eq:decomposition} is a {\em discrete} optimization problem);
usually one has $\dom f_t\subseteq\{0,1\}^{A_t}$.
The arity $|A_t|$ of function $f_t$ can be arbitrarily large, however 
we assume the existence of an efficient {\em min-oracle} that for a given vector $y\in\mathbb R^{A_t}$ computes 
$X\in\argmin\limits_{X\in \dom f_t} \left[f_t(X)+\langle X,y\rangle\right]$
together with the cost $f_t(X)$.

Note, in the case of problem~\eqref{eq:MAPMRF} vector $X$ has the form $X=(X_{va}\::\:v\in [n],a\in D_v)$,
where $X_{va}$ is the indicator variable of the event $[z_v=a]$.
Thus, we have $X\in\{0,1\}^d$ where $d=\sum_{v=1}^n|D_v|$.
Each term $f_t(\cdot)$ is a MAP-MRF problem on a tree-structured graph;
they are chosen in such a way that their sum equals the objective function in~\eqref{eq:MAPMRF}.


\myparagraph{Relaxation of~\eqref{eq:decomposition}}
We form a relaxation of the problem as in~\cite{Swoboda:CVPR19}. 
First, define 
 subsets $\calX_t,\X_t\subseteq \mathbb R^{A_t}\times\mathbb R$ via
$$
\calX_t=\{[X\; f_t(X)]\::\:X\in \dom f_t\}\;\qquad \X_t=\conv(\calX_t)
$$
Throughout the paper we will refer to elements of $\dom f_t$ as ``atoms'' and
to elements as $\calX_t$ as ``extended atoms''.
For a vector $z\in \mathbb R^{A_t}\times\mathbb R$ let $z_\star\in \mathbb R^{A_t}$ be its first $|A_t|$ components 
and $z_\circ\in\mathbb R$ be its last component (so that $z=[z_\star\; z_\circ]$). 
Also denote $\calX=\bigotimes_{t\in  T} \calX_t$ and $\X=\bigotimes_{t\in  T} \X_t=\conv(\calX)$. 
The $t$-th component of vector $x\in\X$
will be denoted as $x^t\in \X_t$.
Problem~\eqref{eq:decomposition} can now be equivalently written as
\begin{equation}
  \label{eq:decomposition1}
  \min_{\substack
{ x \in \X\;,\; X\in\mathbb R^d
\\x^t_\star = X_{A_t}\in\dom f_t\;\; \forall t\in T} 
} \;\;
\sum_{t\in T} x^t_\circ
\end{equation}
The relaxation is formed by removing non-convex constraints $X_{A_t}\in\dom f_t$.
Dualizing constraints $x^t_\star = X_{A_t}$ with Lagrange multipliers $y=(y^t\in\mathbb R^{A_t}\::\:t\in T)$ and eliminating variables $X$ yields Lagrangian
\ifTWOCOLUMNS
	\begin{align*}
	\calL(x;y)&=\sum_{t\in T} (x^t_\circ + \langle x_\star^t,y^t\rangle) \\
	&= \sum_{t\in T}  \langle x^t,[y^t\; 1]\rangle \qquad\quad (x,y)\in\X\times\Y
	\end{align*}
\else
	\begin{equation}
	\calL(x;y)=\sum_{t\in T} (x^t_\circ + \langle x_\star^t,y^t\rangle)= \sum_{t\in T}  \langle x^t,[y^t\; 1]\rangle \qquad\quad (x,y)\in\X\times\Y
	\end{equation}
\fi
where we denoted
\ifTWOCOLUMNS
	\begin{align*}
	&\Y 
	= \left\{ y \:  
	: \:
	\sum\limits_{t\in T_i} y^t_i = 0 \;\;\; \forall i\in[d]\right\}
	\;,\\ &\quad\quad T_i=\{t\in T\::\:i\in A_t\}
	\end{align*}
\else
	\begin{equation}
	\Y 
	= \left\{ y \:  
	: \:
	\sum\limits_{t\in T_i} y^t_i = 0 \;\;\; \forall i\in[d]\right\}
	\;,\qquad T_i=\{t\in T\::\:i\in A_t\}
	\end{equation}
\fi
Note that $\calL(x;y)$ is convex in $x$ and concave in $y$. Define primal and dual objectives via
\begin{equation}
F(x)=\max_{y\in\Y} \calL(x;y)\quad\qquad H(y)=\min_{x\in \X} \calL(x;y)
\end{equation}
Assuming strong duality, the optimal value of the relaxation of~\eqref{eq:decomposition1}
is then equals
\begin{equation}\label{eq:SADDLE}
\min_{x\in\X} F(x)=\max_{y\in \Y} H(y)
\end{equation}
The problem becomes to compute saddle point $(x,y)$ that attains the optimum values in~\eqref{eq:SADDLE}.

\myparagraph{Proximal point algorithm} Next, we review a method for solving~\eqref{eq:SADDLE} proposed in~\cite{Swoboda:CVPR19,KP:ICML21}.
 At each iteration it considers a smoothed version of the problem
with the following Lagrangian and primal and dual objectives:
\ifTWOCOLUMNS
	\begin{align*}
	\calL_{\gamma,\bar y}(x;y)&=\calL(x;y)-\tfrac {||y-\bar y||^2}{2 \gamma}\\
	F_{\gamma,\bar y}(x)&=\max_{y\in\Y} \calL_{\gamma,\bar y}(x;y)\\ 
	H_{\gamma,\bar y}(y)&=\min_{x\in \X} \calL_{\gamma,\bar y}(x;y)=H(y)-\tfrac {||y-\bar y||^2}{2 \gamma}
	\end{align*}
\else
	\begin{equation}
	\calL_{\gamma,\bar y}(x;y)=\calL(x;y)-\tfrac {||y-\bar y||^2}{2 \gamma}
	\end{equation}
	\begin{equation}
	F_{\gamma,\bar y}(x)=\max_{y\in\Y} \calL_{\gamma,\bar y}(x;y)\quad\qquad H_{\gamma,\bar y}(y)=\min_{x\in \X} \calL_{\gamma,\bar y}(x;y)=H(y)-\tfrac {||y-\bar y||^2}{2 \gamma}
	\end{equation}
\fi
Here $\gamma>0$ is a fixed regularization parameter and $\bar y\in\Y$ is the current proximal center.
The algorithm first selects arbitrary vector $y_0\in \Y$, sets $\bar y_0=y_0$,
and then computes iterates $(x_n,y_n,\bar y_n)$ for $n=1,2,\ldots$ using the following equations:
\begin{subequations}
    \begin{eqnarray}
    x_n &\approx& \argmin_{x\in\X} F_{\gamma,\bar y_{n-1}}(x) \label{eq:APPA:a} \\
    y_n &=& \argmax_{y\in\Y} \calL_{\gamma;\bar y_{n-1}}(x;y) \label{eq:APPA:b}  \\
      \bar y_n &=& y_n + \alpha_n (y_n-y_{n-1}) \label{eq:APPA:c} 
    \end{eqnarray} 
\end{subequations}

\cite{Swoboda:CVPR19} used a simple proximal point method (PPA) in which $\alpha_n=0$ for all $n$.
\cite{KP:ICML21} used an accelerated version (APPA) where $\alpha_n$ is given by $\alpha_n=\tfrac{t_n-1}{t_{n+1}}$
and $t_1,t_2,\ldots$ is a sequence satisfying $t_1=1$ and $t_{n-1}^2-t_n^2+t_n>0$ for all $n\ge 2$,
e.g.\ given by the Nesterov's choice $t_{n+1}=(1+\sqrt{1+4t_n^2})/2$.

The main computational bottleneck of the algorithm is the approximate minimization problem in \eqref{eq:APPA:a}.
It is solved using several steps of the Frank-Wolfe (FW) algorithm (or one of its variants) until the FW gap becomes
smaller than some threshold $\varepsilon_n$. (The background on FW is given below).
As shown in~\cite{KP:ICML21}, this guarantees that point $y_n$ satisfies $H_{\gamma,\bar y_{n-1}}(y_n)\ge \max\limits_{y\in\Y} H_{\gamma,\bar y_{n-1}}(y)-\varepsilon_n$, and 
yields accuracy $O(1/n^2)$ after $n$ iterations\footnote{The bound $O(1/n^2)$ holds for the dual, primal and infeasibility gaps, see~\cite{KP:ICML21} for details.}
assuming that $\varepsilon_n=O(1/n^{4+\delta})$ for some $\delta>0$.

\myparagraph{Frank-Wolfe algorithms} For brevity, let us denote $\bar y=\bar y_{n-1}$ and $\tilde F(x)\eqdef F_{\gamma,\bar y}(x)$.
Expanding this, we obtain
\ifTWOCOLUMNS
	\begin{align*}
	&\tilde F(x)=\sum_{t\in T}\left(\tfrac\gamma 2 ||x^t_\star||^2+\langle x^t, [\bar y^t\; 1]\rangle\right) - \sum_{i=1}^d\tfrac{|T_i|}{2\gamma}\nu_i^2,\\
	&\qquad\nu_i=\nu_i(x)=\tfrac{1}{|T_i|}\sum_{t\in T_i}(\gamma\cdot x^t_i + \bar y_i^t)
	\end{align*}
\else
	$$
	\tilde F(x)=\sum_{t\in T}\left(\tfrac\gamma 2 ||x^t_\star||^2+\langle x^t, [\bar y^t\; 1]\rangle\right) - \sum_{i=1}^d\tfrac{|T_i|}{2\gamma}\nu_i^2,\qquad
	\nu_i=\nu_i(x)=\tfrac{1}{|T_i|}\sum_{t\in T_i}(\gamma\cdot x^t_i + \bar y_i^t)
	$$\fi
This is a convex differentiable function over polytope $\X$. Its gradient is given by
$$
\nabla_t \tilde F(x)=[y^t\; 1],\qquad y^t=y^t(x)
=\gamma\cdot \bar x^t_\star + \bar y^t - \nu_{A_t}
$$
Note that $y=y(x)=(y^t(x))_{t\in T}$ is the vector that maximizes $\calL_{\gamma;\bar y_{n-1}}(x;y)$ in eq.~\eqref{eq:APPA:b}.
The basic Frank-Wolfe algorithm with line search~\cite{FWolfe56} minimizes $\tilde F$ by iteratively repeating the following steps: \\
(1) compute $\nabla \tilde F(x)$ at current point $x$; \\
(2) compute $s\in\argmin_{s\in \X}\langle\nabla \tilde F(x),s\rangle$; \\
(3) define $x^\gamma=x+\gamma(s-x)$, compute $\gamma=\argmin\limits_{\gamma\in[0,1]}\tilde F(x^\gamma)$; \\
(4) update $x:=x^\gamma$. \\
Note that step (2) requires calling min-oracles for each subproblem $t\in T$. The quantity
${\tt gap}^{\tt FW}(x;\tilde F)=\langle \nabla \tilde F(x),x-s\rangle=\max_{s\in \X}\langle\nabla \tilde F(x),x-s\rangle$
is called the {\em Frank-Wolfe gap} at $x$. It upper-bounds the suboptimality gap:
${\tt gap}^{\tt FW}(x;\tilde F)\ge \tilde F(x)-\min_{x'\in\X} \tilde F(x')$.

Several techniques have been proposed in the literature for speeding-up the basic algorithm.
We will make use of the following ideas.
\begin{itemize}[topsep=0pt, itemsep=1pt, leftmargin=11pt]
\item Block-Coordinate Frank-Wolfe (BCFW)~\cite{BCFW}. At each step it updates only variables $x^t$ for some fixed $t\in T$
while keeping all other components fixed. More precisely, it computes partial derivatives 
$\nabla_t \tilde F(x)$, then computes $s^t\in\argmin_{s^t\in \X_t}\langle\nabla_t \tilde F(x),s^t\rangle$
and defines $x^\gamma$ via $(x^\gamma)_t=x^t+\gamma(s^t-x^t)$ and $(x^\gamma)_{t'}=x^{t'}$ for $t'\ne t$. The rest is as above.
\item Caching atoms~\cite{Joachims2009,MP-BCFW,Osokin:ICML16}. The idea is store ``atoms'' $s^t\in\calX_t$ returned by the $t$-th min-oracle
in cache  $\calW_t\subset\calX_t$.
Iterations are then divided into ``exact'' iterations which call the (expensive) min-oracle that optimize over $\X_t$,
and ``approximate'' iterations which optimize over $\calW_t$.
\item Optimizing the objective over extended atoms in $\calW_t$
(i.e.\ subject to the constraint $x^t\in{\tt conv}(\calW_t)$, assuming that current $x$ satisfies this constraint).
This amounts to minimizing a smooth convex function (quadratic in our case) over a simplex.
An example is the BCG method~\cite{BCG}.
\item In-face Frank-Wolfe directions~\cite{FW:in-face}. In addition to ``regular'' FW steps, this method performs the following
operations: find a face $\X'\subseteq\X$ of polytope $\X$ containing current point $x$ (e.g.\ the minimal such face), and run several FW steps to minimize function $\tilde F$
over $\X'$. This may lead to a speed-up if oracle calls over $\X'$ are faster than over $\X$.
We will say that $\X'$ is obtained from $\X$ via ``contraction'', and refer to FW steps over $\X'$ as ``contracted steps''.
The same terminology will be used for subproblems $\X_t$.
\end{itemize}


\section{Our implementation}

Following~\cite{Joachims2009,MP-BCFW,Osokin:ICML16},
 for each 
subproblem $t$ we maintain cache $\calW_t\subseteq \calX_t$ of extended atoms.
We  implemented two versions: \\
(1) $x^t$ is represented as a convex combination of elements in $\calW_t$.
Then we store coefficients of this combination. Extended atoms with zero weight are immediately removed. \\
(2) Cache size is limited  by a constant (10 in our implementation).
As in previous works, we maintain a timestamp for each atom, updating it whenever the atom is ``accessed''
(i.e.\ returned as optimal in one of the operations). If the cache becomes too big, we remove an atom with the oldest timestamp. \\
Let ``${\tt conv}$'' be the flag that specifies the version: 
${\tt conv}=\true$ corresponds to option (1) and
${\tt conv}=\false$ corresponds to option (2).

\myparagraph{In-face FW directions}
As stated in the introduction, the main motivation of this paper is to incorporate in-face FW directions
for combinatorial subproblems into the framework of~\cite{Swoboda:CVPR19,KP:ICML21}.
We focus on faces of polytopes $\X'_t\subseteq\X_t$ that are specified by constraints of the form $x^t_i={\tt const}_i$ for some  $i\in A_t$.~\footnote{Our implementation
also supports constraints of the form $x^t_i=x^t_j$ for some $\{i,j\}\subseteq A_t$, but we haven't used it in experiments.}
The following issues should be taken into account:
\begin{itemize}[topsep=0pt, itemsep=1pt, leftmargin=11pt]
\item Fixing some variables may split subproblems into several independent subproblems (that we call ``contracted subproblems''),
and so variables of $x^t$ decouple into independent blocks. 
Optimizing over $x^t\in\X'_t$ should be done via the block-coordinate version of FW.
\item Atoms in the cache $\calW_t$ for the original subproblem $t$ should be transformed to atoms of contracted subproblems,
and vice versa. Ideally, the time for processing each atom should depend on the size of the contracted subproblems (in the case when $|\calW_t|\gg 1$).
\end{itemize}
To address these issues,  we will describe an abstract data structure called {\tt Subproblem}.
We will later show to implement it for subproblems $t$
corresponding to the MAP-MRF inference problems on tree-structured graphs.

We distinguish between ``parent'' and ``child'' subproblems. The former are added to the solver during initialization.
The user should specify mapping $A_t$ when adding parent subproblem (which is an array  size $|A_t|$). 
Now suppose that the algorithm decides to contract parent subproblem $t$.
Given current vector $x$,
the user should first partition set $A_t$ as
$
A_t=A_{t0} \cup A_{t1} \cup \ldots \cup A_{tk}
$ for some $k\ge 1$ where $A_{t0}$ is the set of components that will be fixed to their current values,
and $A_{t1},\ldots,A_{tk}$ correspond to the independent subproblems.
Let us write vectors $z\in\mathbb R^{A_t}$ as $(z_{A_{t0}},z_{A_{t1}},\ldots,z_{A_{tk}})$.
We require the contraction operation to satisfy the following conditions: \\
(1)
There should exist functions $f_{ti}:\{0,1\}^{A_{ti}}\rightarrow\mathbb R\cup\{+\infty\}$
so that
\ifTWOCOLUMNS
	\begin{align*}
	&f_t(x_{A_{t0}},X^1,\ldots,X^k)=f_{t1}(X^1)+\ldots+f_{tk}(X^k)\\
	&\hspace{50pt}\forall (X^1,\ldots,X^k)\in\dom f_{t1}\times\ldots\times\dom f_{tk}
	\end{align*}
\else
	$$
	f_t(x_{A_{t0}},X^1,\ldots,X^k)=f_{t1}(X^1)+\ldots+f_{tk}(X^k)\qquad\quad\forall (X^1,\ldots,X^k)\in\dom f_{t1}\times\ldots\times\dom f_{tk}
	$$
\fi
(2) For each $i\in[k]$ there should exist vector $\bar x^i\in\X_{ti}$ with $\bar x^i_\star=x^t_{A_{ti}}$
so that $\bar x^1_\circ+\ldots+\bar x^k_\circ=x_\circ$,
where we denoted
$
\calX_{ti}=\{[X\; f_{ti}(X)]\::\:X\in \dom f_{ti}\}$
and
$\X_{ti}=\conv(\calX_{ti})
$. \\
(3) Extended atom $s^t\in\argmin_{s^t\in \calX_t}\langle\nabla_t \tilde F(x),s^t\rangle$
should satisfy $s^t_{A_{ti}}\in\dom f_{ti}$ for each $i\in[k]$.

Accordingly, each parent subproblem should implement function ${\tt Contract}(x^t,s^t)$.
This function should find a decomposition into $k$ child subproblems as above,
and return these subproblems together with mappings $A_{t1},\ldots,A_{tk}$
and values $\bar x^1_\circ,\ldots,\bar x^k_\circ$.
Note that the latter values are actually not used if ${\tt conv}=\true$,
since in this case the solver has enough information to recompute them. 
However, if ${\tt conv}=\false$ then the solver does need values $\bar x^1_\circ,\ldots,\bar x^k_\circ$.

Below we will use index symbol $t$ to denote both parent and child subproblems (together with notation $A_t$, $x^t$, $f_t$, etc).
Note, for child subproblems we have $t=t'i$ for some parent subproblem~$t'$.

\myparagraph{Compact representation of atoms}
In many combinatorial problems, atoms (i.e.\ elements of $\dom f_t$) 
can be described in a compact way. Accordingly, for each subproblem we introduce type {\tt Atom} whose implementation
should be specified by the user. 
The user should implement the following basic functions:
\begin{itemize}[topsep=0pt, itemsep=1pt, leftmargin=11pt]
\item ${\tt AtomToVector}(a)$: compute extended atom $[a\; f(a)]\in\calX_t$ corresponding to atom $a$.
\item ${\tt DotProduct}(a,g)$: compute inner product $\langle a,g \rangle$ of atom $a$ and   vector $g\in\mathbb R^{A_t}$.
\item ${\tt WeightedDotProduct}(a,b)$: compute $\sum_{i\in A_t} (1-\tfrac 1{|T_i|})a_ib_i$ for atoms $a,b$.
\item ${\tt MinOracle}(g)$: return atom $a\in\argmin_{a} (\langle a,g\rangle + f_t(a))$.
\end{itemize}
Next, we discuss how to transform atoms during contractions.
Each parent subproblem $t$ must store a current atom denoted as $a^t$,
and implement functions $t\!::\!{\tt SetAtom}(a)$ and $t\!::\!{\tt GetAtom}()$
that respectively set and return $a^t$.
Child subproblems $ti$ must also implement these two functions, but they are now defined as follows.
$ti\!::\!{\tt SetAtom}(a)$ should set components of the parent atom $a^t$ to~$a$, i.e.\ update
$a^t_{A_{ti}}:=a$.
Similarly, $ti\!::\!{\tt GetAtom}()$ should return atom $a^t_{A_{ti}}$ of the child subproblem
(or {\tt NULL}, if $a^t_{A_{ti}}\notin\dom f_{ti}$).

We now describe how  these functions are used.
Consider subproblem $t$ which has been contracted to child subproblems $t1,\ldots,tk$
via the call ${\tt Contract}(x,s)$ (for brevity, we write $x,s$ instead of $x^t,s^t$).
By construction, atom $a^t$ will satisfy $a^t_{A_{t0}}=s_{A_{t0}}$ where $A_{t0}$ are the variables that have been fixed.
Now suppose that the solver decides to contract $t$ again (i.e.\ recompute child subproblems).
The solver first computes atom $s'\in\argmin_{s'\in \calX_t}\langle\nabla_t \tilde F(x'),s'\rangle$
where $x'\in\X_t$ is the current vector for $t$.
Then it calls ${\tt Contract}(x',s')$
which returns new child subproblems $t1',\ldots,tk'$.
If ${\tt conv}=\false$ then caches $\calW_{t1'},\ldots,\calW_{tk'}$ are set by repeating the following steps: \\
(1) for each $i\in[k]$ pick atom $a^i\in\calW_{ti}$ (in a round-robin fashion), call $ti\!::\!{\tt SetAtom}(a^i)$; \\
(2) for each $i'\in[k']$ call $a^{i'}\leftarrow ti'\!::\!{\tt GetAtom}()$,
add $[a^{i'}\; f_{ti'}(a^{i'})]$  to $\calW_{ti'}$ (if $a^{i'}\ne{\tt NULL})$.

\noindent
If ${\tt conv}=\true$
then we  use $ti\!::\!{\tt SetAtom}(\cdot)$ and $ti'\!::\!{\tt GetAtom}()$ in a similar way
to restore the desired invariant, i.e.\ make sure that $x^{ti'}$ is a convex combination of current atoms for each $i'\in[k']$;
details are omitted.
In the end of the contraction operation we call $t\!::\!{\tt SetAtom}(s')$.

\subsection{MAP-MRF tree subproblems}\label{sec:treeMAP}
We now detail the implementation of {\tt Subproblem} in the case when $f_t(\cdot)$ corresponds to the MAP-MRF inference problem on a tree-structured graph.
In other words, we assume that $f_t(\cdot)$ corresponds to problem~\eqref{eq:MAPMRF} (with labelings $z$  represented via vectors $X\in\{0,1\}^d$) 
in which the graph $([n],\calE)$ is a tree. Recall than $d=\sum_{v=1}^n|D_v|$.


An atom is stored in a natural way as a vector of size $n$ (which can be much smaller than ${\tt dim}(X)=d$).
Next, we discuss the implementation of function ${\tt Contract}(x,s)$ for vector $x\in\mathbb R^{d+1}$ and atom $s\in D_1\times\ldots \times D_n$.
First, we find pairs $(v,a)$ satisfying $x_{va}=0$ and $s_v\ne a$, and force constraint $X_{va}=0$ for such pairs
(assuming that there are least 25\% of such pairs; otherwise we do not contract).
Let $\calV_{\tt fixed}\subseteq[n]$ be the set of nodes $v$ for which $|D_v|-1$ variables $X_{va}$ has been forced to $0$,
and let $\calV=[n]-\calV_{\tt fixed}$. Clearly, the problem splits into $k\ge 0$ independent subproblems
corresponding to the trees $\calG_1,\ldots,\calG_k$ of the induced forest $\calG[\calV]$.
Each node of the child subproblem has a pointer to the corresponding node of the parent subproblem (and vice versa).
Clearly, this allows an efficient implementation of functions ${\tt GetAtom}()$ and ${\tt SetAtom}(\cdot)$.

As discussed in the previous section, we need to be able to compute values $\bar x^k_\circ$ for each child subproblem $i\in[k]$
(in the case when ${\tt conv}=\false$).
For that we need to know the cost of (fractional) vector $x$ restricted to tree~$\calG_i$.
This can be easily computed if we know the ``fractional cost'' $\xi_{uv}$ of each edge $uv\in\calE$ in $x$.
This fractional cost must have the following form for some vector $(\xi_{ua;vb})_{a,b}$: 
\begin{subequations}\label{eq:xiuv}
\begin{align}
&\xi_{uv} = \sum_{a\in D_u,b\in D_v} f_{uv}(a,b)\xi_{ua;vb} \\
& \sum_{b\in D_v} \xi_{ua;vb} = x_{ua} \qquad \forall a\in D_u \\
& \sum_{a\in D_u} \xi_{ua;vb} = x_{vb}  \qquad \forall b\in D_v 
\end{align}
\end{subequations}
Accordingly, we additionally maintain $|\calE|$ real numbers for the current vector $x$ as described below.

One possibility would be to increase the dimensions of vector $X$,
i.e.\ let $X=(X_{va}\::\:v\in [n],a\in D_v)\sqcup (X_{uv}\::\:uv\in\calE)$.
If $X$ corresponds to atom $z=(z_1,\ldots,z_n)$ then $X_{uv}=f_{uv}(z_u,z_v)$.
The solver would then automatically maintain the desired fractional cost $\xi_{uv}=X_{uv}$ for each edge $uv$.
We opted for another approach: we store costs $\xi_{uv}$ internally at the given subproblem.
Note that solver changes current $x$ via updates of the form $x:=\alpha_0 x + \sum_{i=1}^r \alpha_i{a^i}$
where $a^1,\ldots,a^r$ are extended atoms and $\alpha_0,\alpha_1,\ldots,\alpha_r$ are non-negative coefficients that sum to 1.
Before applying such update, the solver calls function ${\tt Update\_\,x}(\alpha_0,x,\alpha_1,a^1,\ldots,\alpha_r,a^r)$ for the the given subproblem,
and this function updates $\xi_{uv}$ for all edges $uv$. This allows the solver to manipulate with vectors of smaller sizes.

\myparagraph{Optimal transportation problem}
Note that minimizing the value of $\xi_{uv}$ over $(\xi_{ua;vb})_{a,b}$ in~\eqref{eq:xiuv}
(for fixed vectors $(x_{ua})_a$ and $(x_{vb})_b$) is a classical optimal transportation (OT) problem.
It is usually small-scale; we solve it via a successive shortest path algorithm. We use in two ways.

First, we observed that sometimes the FW approach with option ${\tt conv}=\false$
makes some components of vector $x$  extremely small, e.g.\ $10^{-20}$ or $10^{-30}$.
This can be explained as follows. Suppose that the value $x_{va}$ becomes positive in the early stage of the algorithm,
but later on all atoms $z$ satisfy $z_v\ne a$. Furthermore, suppose that optimizing the objective over ${\tt conv}(x^t\cup\calW_t)$
always assigns a positive weight to $x^t$, and this weight is smaller than some constant $c<1$ sufficiently often.
Then after $n$ such steps $x_{va}$ becomes smaller than $c^n$.

To address this issue, we do the following from time to time (if ${\tt conv}=\false$):
we go through all subproblems $t$ and update $x^t$ by calling $x^t\leftarrow t\!::\!{\tt Adjust\_x}(x^t)$.
The latter function is implemented as follows for tree-structured MAP-MRF subproblems:
we change all fractional components smaller than $10^{-8}$ to zero, renormalize unary fractional components for each node
so that they sum to 1, and then compute the fractional cost $\xi_{uv}$ for each edge $uv$ and the total cost by solving OT problems.
The frequency of such updates is set so that their time is at most 20\% of the total time.
More precisely, we run the ${\tt Adjust\_x}$ loop after each update of proximal center $\bar y$
assuming that $\tau > 5\tau_0$ where $\tau_0$ is the time spent in the previous ${\tt Adjust\_x}$ loop
and $\tau$ is the time that elapsed after the last ${\tt Adjust\_x}$ loop.
This scheme improved the performance of in-face FW directions for some instances, since it led to smaller contracted subproblems.
For the algorithm without in-face FW directions this scheme affected the performance only marginally.

Following~\cite{SavchynskyySchmidt:13}, we also use the OT procedure for computing feasible primal solutions of relaxation~\eqref{eq:BLP}.
Note that vectors $x\in\X$ do not directly give such a solution,
since for each node $v\in [n]$ and label $a\in D_v$
variables $x^t_{va}$ will in general be different for different subproblems $t$ containing $v$.
To circumvent this problem, we compute the average of values $x^t_{va}$ over $t$ containing $v$;
this gives a feasible fractional solution for nodes. The corresponding fractional solution for edges (and the LP cost of this solution)
is then computed by solving OT problems.

\subsection{Perfect Matching subproblems}
We also experimented with subproblems corresponding to functions of the following form: $f_t(X)=\langle C,X\rangle$
if $X\in\{0,1\}^\calE$ is a perfect matching in an undirected graph $(\calV,\calE)$, and $f_t(X)=+\infty$ otherwise.
We obtained such subproblems by considering quadratic pseudo-boolean optimization problems 
(minimize $f(z)=\sum_{v} c_v z_v + \sum_{uv} c_{uv} z_u z_v$ 
over $z\in\{0,1\}^n$),
and tightening it with planar subproblems~\cite{Yarkony:UAI11}; the latter are solved by a reduction to the perfect matching problem.
A face of the perfect matching polytope can be described by constraints of the form 
\begin{subequations}
\begin{eqnarray}
X(\calE_0)&=&{\bf 0} \\
X(\delta S)&=&1\qquad\forall S\in\calS \label{eq:PM-S} 
\end{eqnarray}
\end{subequations} 
where $\calE_0\subseteq\calE$, $\calS$ is a laminar family of odd-cardinality subsets of $\calV$, and $\delta S$ is the set of edges of the cut $(S,\calV-S)$.
Given a fractional point~$X$ in the polytope, computing the minimal face containing~$X$ boils down to computing a Gomory-Hu tree in graph $(\calV,\calE)$ with appropriate weights~\cite{PadbergRao82}.

In our informal experiments (not reported here) we were not able to obtain a speed-up with in-face FW directions. First, we observed that set $\calE_0$ is usually empty (or very small).
This could be due to the way that our planar subproblems are formed (we add to planar subproblems  odd cycles that are currently violated).
The runtime of Gomory-Hu computations was negligible using the recent code of~\cite{OC,OC:experimental}.
However, adding constraints~\eqref{eq:PM-S} did not not make the Blossom V code~\cite{blossom5} faster.
Note that we incorporated constraints~\eqref{eq:PM-S}  by adding a large constant to edge weights of edges in sets $\delta S$ for $S\in\calS$.

We conjecture that in-face FW directions could still be beneficial
 either (i) in applications where set $\calE_0$ is large, or (ii) if the perfect matching
 code is modified so that constraints~\eqref{eq:PM-S} are incorporated more directly. We leave this as a future work.

\subsection{Frank-Wolfe algorithm} 

In this section  we describe the FW version that we have implemented.
One basic operation is to minimize $\tilde F(\ldots, x^t,\ldots)$  over $x^t\in{\tt conv}(Z)$
for some finite set $Z=\{z^1,\ldots,z^k\}$;
all other components except for $x^t$ are fixed.\footnote{If
${\tt conv}=\true$ then $Z=\calW_t$, and if ${\tt conv}=\false$ then $Z=\calW_t\cup\{x^t\}$ where $x^t$ is the current vector.}
Denote $x^t=\sum_{i=1}^k \alpha_i z^i$, then equivalently we need to minimize quadratic function $f(\alpha)=\tfrac 12 \alpha^T A \alpha + b^T \alpha$
over simplex $\alpha\in\{\alpha\in\mathbb R^k_{\ge 0}\::\:\sum_{i=1}^k\alpha_i=1\}$.
Define vector $\lambda\in\mathbb R^{A_t}$ via
$
\lambda_i=\sum_{t\in T_i}\tfrac {\gamma x^t_i}{|T_i|} - \nu_i(x)
$ where $x$ is the current vector, then  matrix $A\in \mathbb R^{k\times k}$ and vector $b\in\mathbb R^k$ are given by
$
A_{ij}=\gamma\cdot {\tt WeightedDotProduct}(z^i,z^j)
$
and
$
b_i=\langle z^i,[\bar y^t\; 1]\rangle + \langle z^i_\ast, \lambda \rangle
$.
Since values ${\tt WeightedDotProduct}(a,b)$ for $a,b\in\calW_t$  are needed multiple times,
we store them for each cache $\calW_t$ in a matrix of size $|\calW_t|\times |\calW_t|$.

To minimize  $f(\cdot)$, we implemented a projected conjugate gradient descent (PCG) method
with the following modification. Suppose that a step takes the procedure outside the simplex, i.e.\ the weight $\alpha_i$ of
point $z^i$ becomes negative. We then stop at the boundary (making $\alpha_i=0$), remove $z^i$ from $Z$, and restart PCG. 
Thus, each step either makes $Z$ smaller (which can happen at most $|Z|-1$ times)
or makes a sufficient progress on the objective (as guaranteed by PCG).
We stop when the FW gap is reduced by a factor of 10.
This can be viewed as a particular implementation of the {\em Simplex Descent Oracle} in the BCG method~\cite{BCG}.


We use three types of operations: {\tt CacheLoop}, {\tt OracleLoop} and {\tt ContractLoop}. In the first two we cycle through child subproblems in a random
order, and in the third one we cycle through parent subproblems in a random order.
For each subproblem $t$ we do the following.
\begin{sloppypar}
\begin{itemize}[topsep=0pt, itemsep=1pt, leftmargin=11pt]
\item {\tt CacheLoop}: run Simplex Descent Oracle ({\tt SiDO}) described above (for $Z=\calW_t$ or $Z=\calW_t\cup\{x^t\}$, depending on flag ${\tt conv}$).  
\item {\tt OracleLoop}: run {\tt SiDO}; compute $\nabla_t \tilde F(x)$; call $t$-th oracle: $s^t\leftarrow\argmin_{s^t\in \X_t}\langle\nabla_t \tilde F(x),s^t\rangle$;
add $s^t$ to $\calW_t$; run {\tt SiDO} again.
Note that {\tt SiDO} is usually much faster than the min-oracle (especially if $|\calW_t|$ is small), so using
it before the oracle call is a cheap way to get a better gradient $\nabla_t \tilde F(x)$.
\item {\tt ContractLoop}: compute $\nabla_t \tilde F(x)$; call $t$-th oracle: $s^t\leftarrow\argmin_{s^t\in \X_t}\langle\nabla_t \tilde F(x),s^t\rangle$;
call ${\tt Contract}(x^t,s^t)$; update child subproblems and transform atoms as described in the previous section.
\end{itemize}
\end{sloppypar}

\begin{figure*}[ht]
{\small \vspace{-50pt}
\begin{center}

		\begin{tabular}{@{\hspace{-22pt}}cc@{\hspace{-17pt}}c @{\hspace{-17pt}}c@{\hspace{-17pt}}c} 
		\raisebox{60pt}{\begin{tabular}{l} X-axis: runtime (in seconds) \\ Y-axis below 0: lower bound \\ Y axis above 0: upper bound \end{tabular}} &
			\raisebox{20pt}{\includegraphics[scale=0.6,trim=61 180 250 44,clip=true]{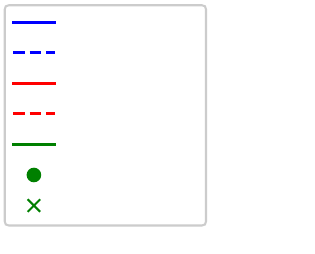} 
			\begin{picture}(0, 0)
				\put(-70,65){\scriptsize ${\tt FW}^\ast_{\tt conv}$} 
				\put(-70,56){\scriptsize ${\tt FW}_{\tt conv}$} 
				\put(-70,47){\scriptsize ${\tt FW}^\ast$} 
				\put(-70,38){\scriptsize ${\tt FW}$} 
				\put(-70,29){\scriptsize TRW-S} 
				\put(-70,20.5){\scriptsize ADSAL} 
				\put(-70,12){\scriptsize AD3} 
				\put(-75,79){\footnotesize legend} 
					\end{picture}
			} &
			\includegraphics[scale=0.32]{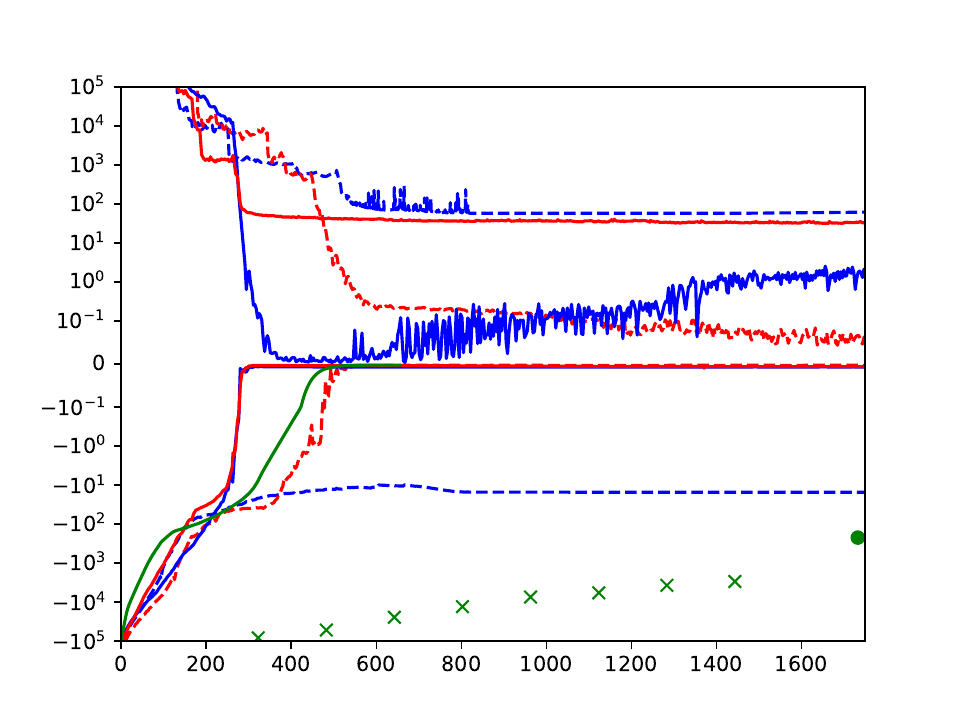} \begin{picture}(0, 0) \put(-85,90){\footnotesize family-gm} \end{picture} &
			\includegraphics[scale=0.32]{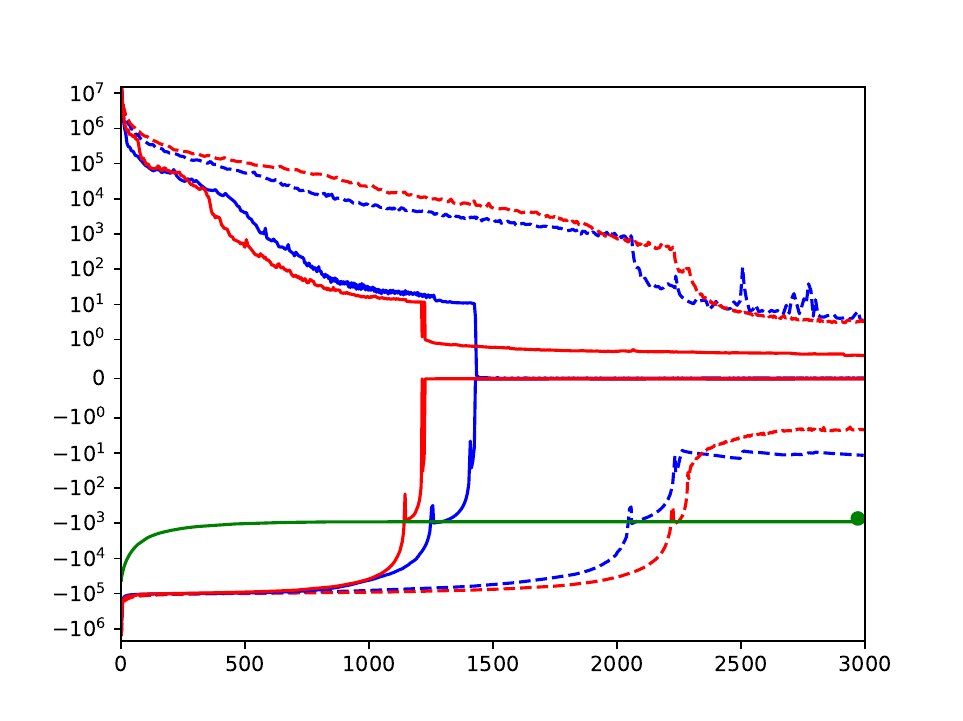} \begin{picture}(0, 0) \put(-85,90){\footnotesize pano-gm} \end{picture} 
		\end{tabular}\vspace{-5pt}

		\begin{tabular}{@{\hspace{-22pt}}c@{\hspace{-17pt}}c @{\hspace{-17pt}}c@{\hspace{-17pt}}c} 
			\includegraphics[scale=0.32]{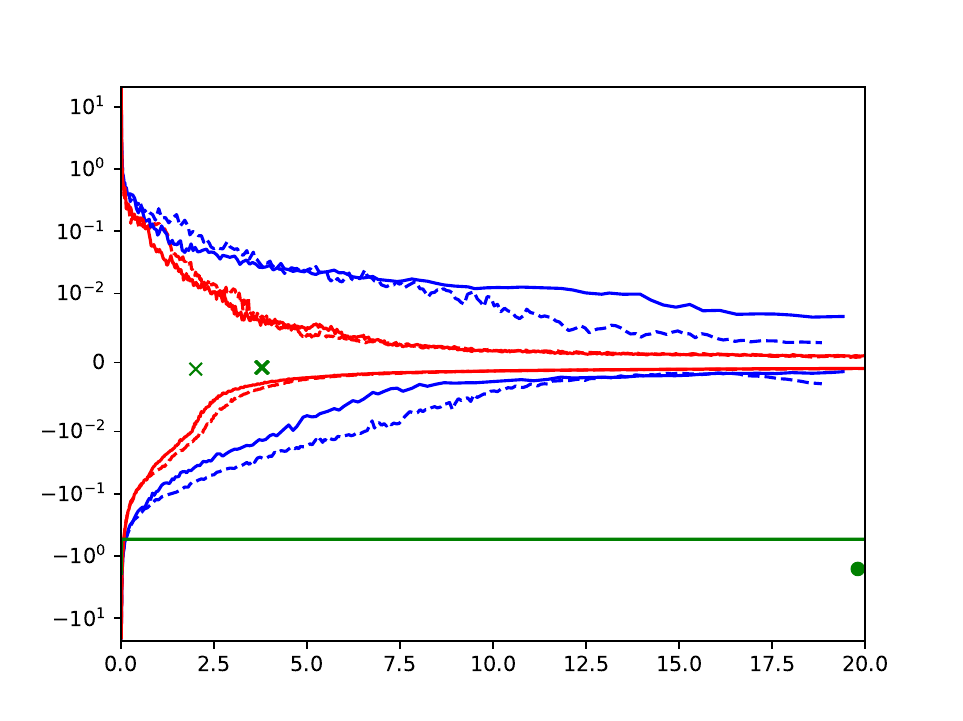} \begin{picture}(0, 0) \put(-85,90){\footnotesize matching0} \end{picture} &
			\includegraphics[scale=0.32]{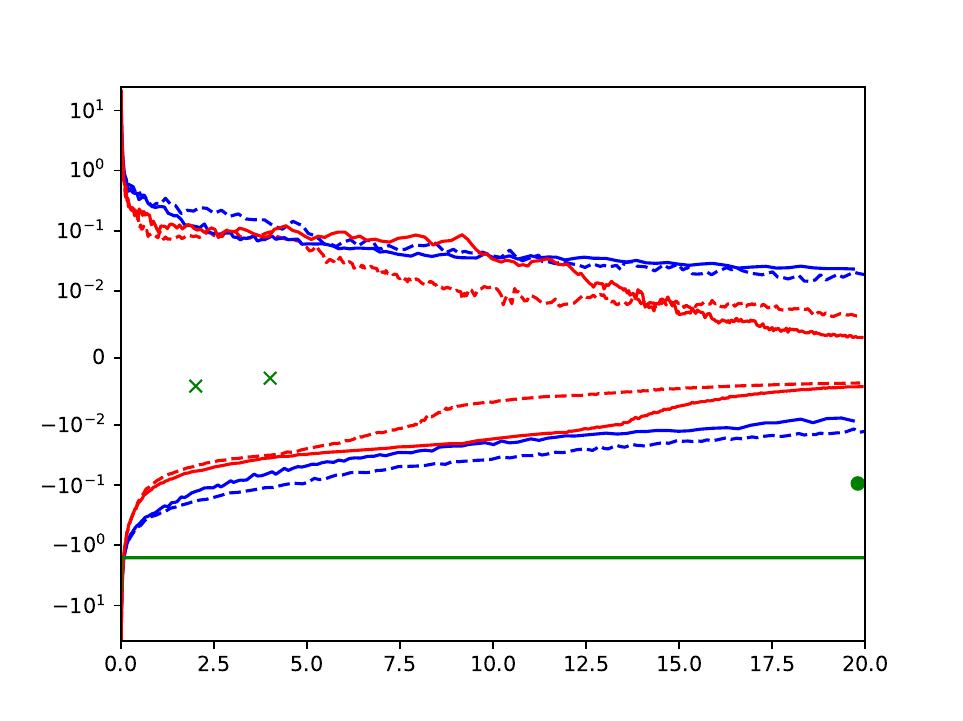} \begin{picture}(0, 0) \put(-85,90){\footnotesize matching1} \end{picture} &
			\includegraphics[scale=0.32]{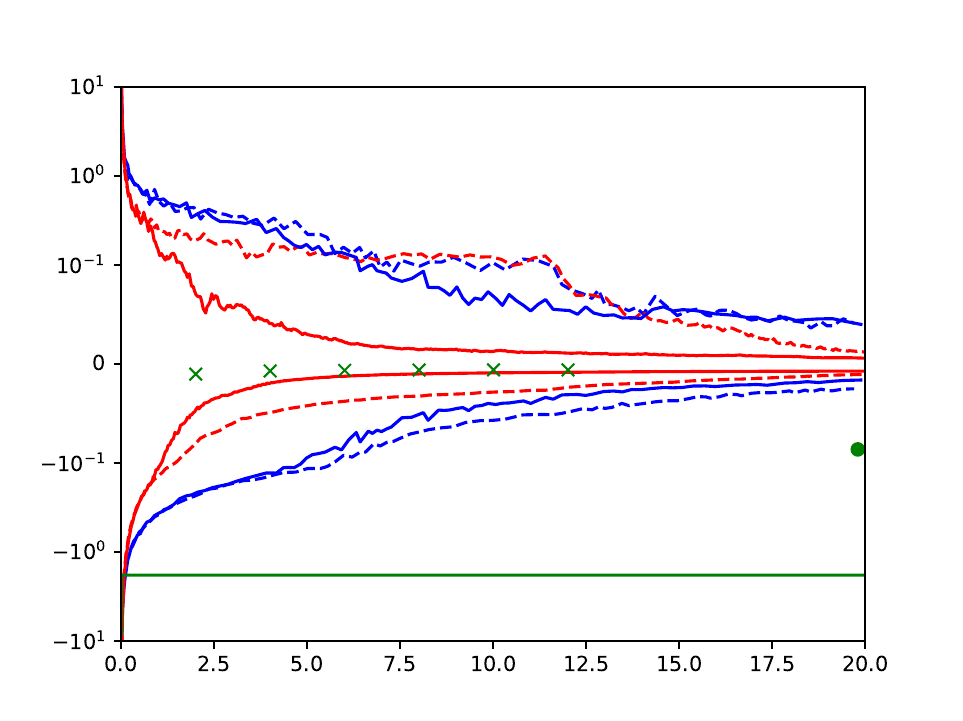} \begin{picture}(0, 0) \put(-85,90){\footnotesize matching2} \end{picture} &
			\includegraphics[scale=0.32]{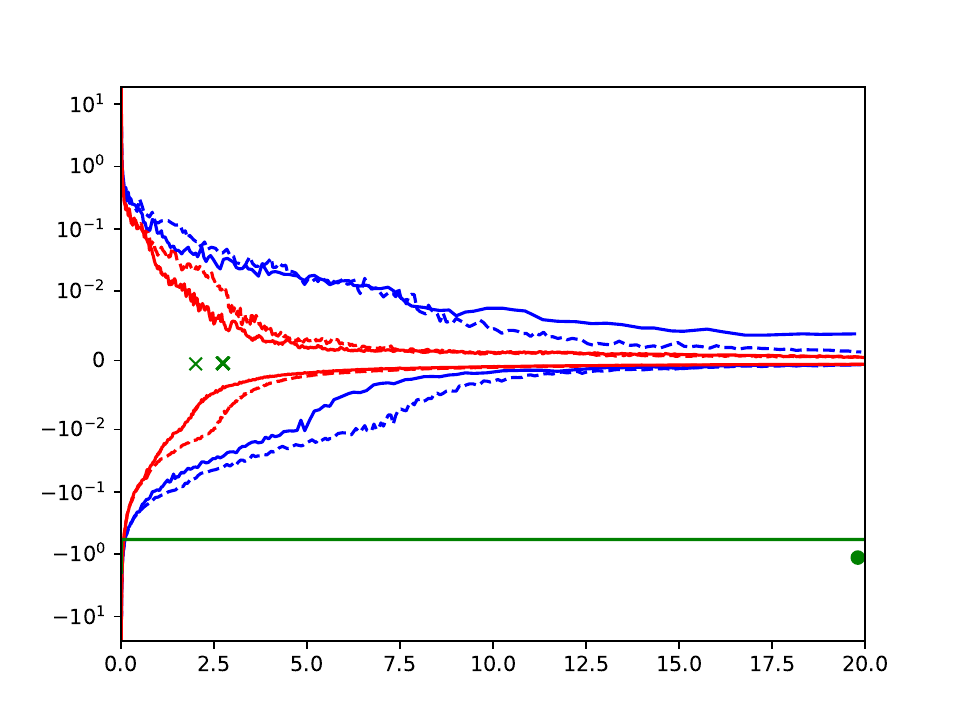} \begin{picture}(0, 0) \put(-85,90){\footnotesize matching3} \end{picture} 
		\end{tabular}\vspace{-5pt}
		
		\begin{tabular}{@{\hspace{-22pt}}c@{\hspace{-17pt}}c @{\hspace{-17pt}}c@{\hspace{-17pt}}c} 
			\includegraphics[scale=0.32]{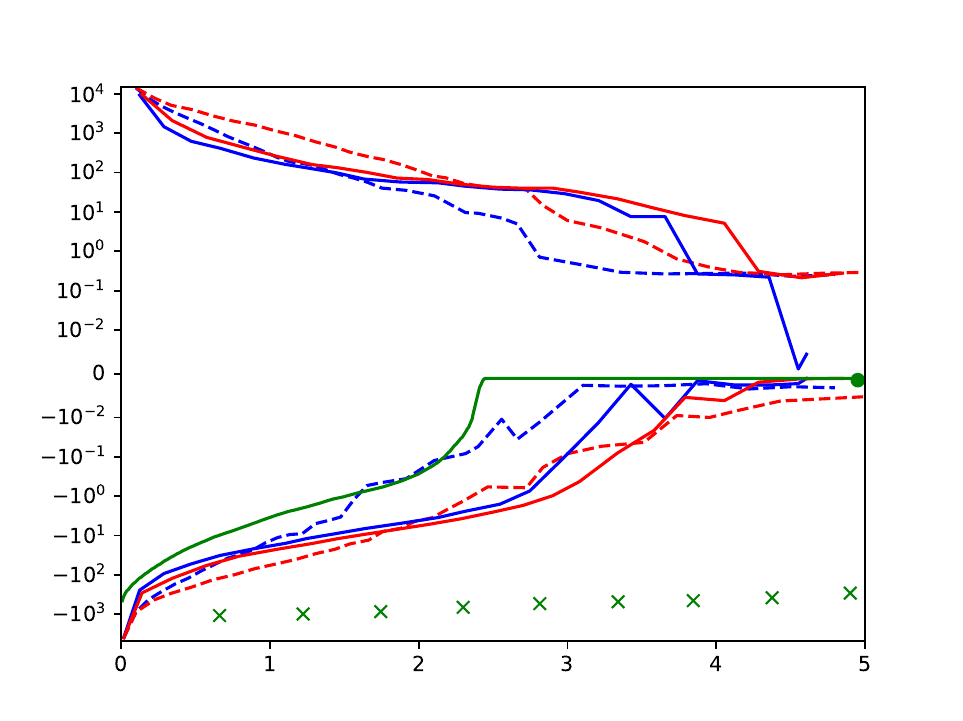} \begin{picture}(0, 0) \put(-85,90){\footnotesize objseg-349} \end{picture} &
			\includegraphics[scale=0.32]{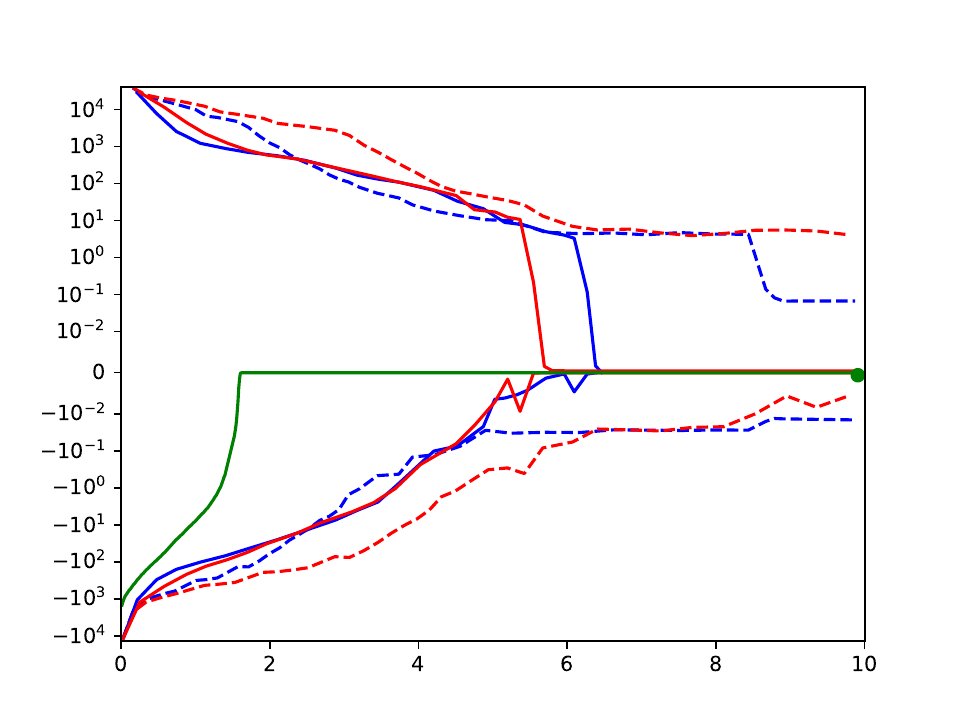} \begin{picture}(0, 0) \put(-85,90){\footnotesize objseg-353} \end{picture} &
			\includegraphics[scale=0.32]{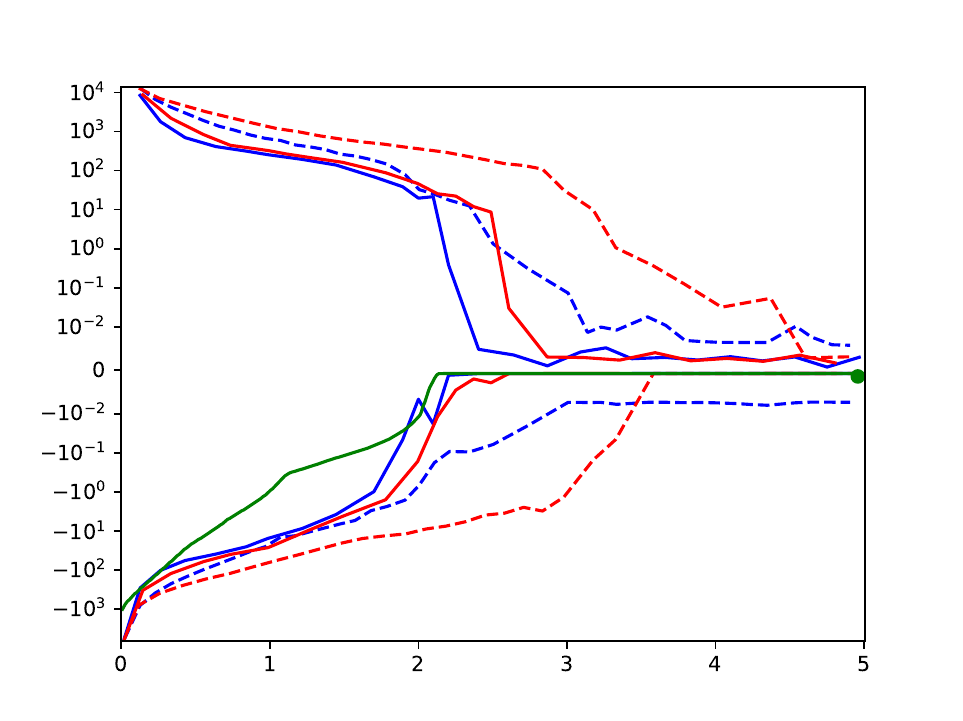} \begin{picture}(0, 0) \put(-85,90){\footnotesize objseg-358} \end{picture} &
			\includegraphics[scale=0.32]{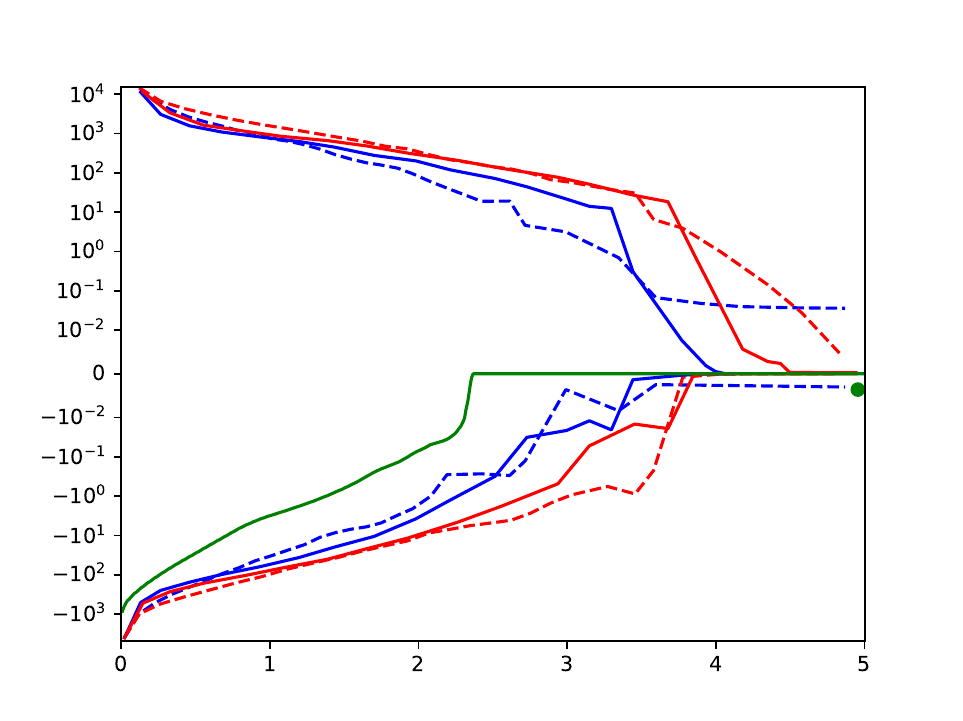} \begin{picture}(0, 0) \put(-85,90){\footnotesize objseg-416} \end{picture} 
		\end{tabular}\vspace{-15pt}

	\end{center}
	}
	\caption{Top row: `mrf-photomontage'  ($\gamma=10$). Middle row: `matching' ($\gamma=1$). Bottom row: `object-seg' ($\gamma=1$).
        }
\label{fig:plots}
\end{figure*}

An important question is how often each operation should be run. In general, {\tt CacheLoop} makes a smaller progress
compared to {\tt OracleLoop} (i.e.\ improves the objective by a smaller amount),
and similarly  in-face FW steps can be expected to make a smaller progress compared to regular steps.
However, the former operations can be much faster. We thus argue that the criterion for choosing the operation should
be based on the actual runtimes.
To choose between {\tt CacheLoop} and {\tt OracleLoop}, we use the same method 
 as in~\cite{MP-BCFW}. Namely, we define procedure {\tt InnerIteration} as follows.
It calls  {\tt OracleLoop} once and then {\tt CacheLoop} one or more time.
Let $\Delta_i$ be the improvement in the objective after the $i$-th call, and $\tau_i$ be the time after the $i$-th call
(both relative to the beginning of {\tt InnerIteration}). We terminate the iteration if after the $i$-th call we have
$
\frac{\Delta_i}{\tau_i}\le \frac{\Delta_{i-1}}{\tau_{i-1}}
$ (for $i\ge 2$). We use a similar technique to choose between {\tt InnerIteration} and {\tt ContractLoop}:
we define procedure {\tt OuterIteration} that calls {\tt ContractLoop} once and then {\tt InnerIteration} one or more time,
using a similar termination criterion.

\myparagraph{Proximal Point Algorithm}
Following~\cite{KP:ICML21}, we use inexact Accelerated Proximal Point Algorithm (APPA). The $n$-subproblem is solved to accuracy $\varepsilon_n=\varepsilon_0/n^2$
where $\varepsilon_0$ is the initial FW gap. Unlike~\cite{KP:ICML21}, we restart the APPA method if the the primal objective $\tilde F(x)$ becomes worse;
this is a well-known strategy in the case of (exact) accelerated algorithms~\cite{AdaptiveRestarts}.

\myparagraph{Numerical stability} As we described in Section~\ref{sec:treeMAP},
if ${\tt conv}=\false$ then some components of vector $x$ may become extremely small.
A similar issue may happen with ${\tt conv}=\true$: weights of some atoms may become small, e.g.\ $1e^{-20}$.
As a result, the algorithm may get stuck, or $|\calW_t|$ may grow very large.
To avoid this issue, we use a similar technique as in Section~\ref{sec:treeMAP}:
we simply remove atoms whose weight is smaller than $1e^{-8}$.

\section{Experimental results}\label{sec:experiments}

In this section we test different algorithms on (pairwise) MAP-MRF inference problems 
from the OPENGM2 benchmark study~\cite{OpenGMBenchmark}.
Unlike~\cite{OpenGMBenchmark}, we do not aim to obtain the best possible solution of~\eqref{eq:MAPMRF}.
Instead, we investigate what is the best algorithm for solving its BLP relaxation~\eqref{eq:BLP}.
Note that this is a well-defined subproblem, and used as a subroutine for other solvers in~\cite{OpenGMBenchmark}
(branch-and-bound techniques, techniques combining LP solvers with persistency criteria, etc).

Accordingly, we included TRW-S method~\cite{Kolmogorov2006} and ADSAL method~\cite{AdaptiveDiminishingSmoothingSavchynskyyUAI}
that seem to dominate other LP solvers according to the results of~\cite{OpenGMBenchmark}.
On a subset of instances we also ran the Gurobi solver (ver. 10.0.0) and the ADMM-based ``AD3'' method~\cite{Martins:ICML11}.\footnote{
The code of \url{https://github.com/andre-martins/AD3} includes parameter
$\eta$ and also a possibility to adapt this parameter during optimization. It appears that \cite{OpenGMBenchmark}
used the default value of $\eta$, namely $\eta=0.1$. We noticed that $\eta$ was changed very rarely during optimization
(usually halved once). We turned off {\tt adapt\_eta} flag and hand-picked $\eta$ which gave the best result.

Note, objseg-349 uses Potts interaction potentials. This is exploited in our implementation (via distance transforms)
but not in AD3 and not in the reduction to Gurobi.
}
Gurobi runtimes were as follows (in seconds, using 4 physical cores; for the first instance Gurobi  aborted, probably because of requesting too much memory:
we used a machine with 16Gb  RAM, out of which 12Gb was free):

\begin{tabular}{|c|c|c|c|}\hline
family-gm & matching0 & objseg-349 & 1CKK  \\
- & 1.36 & 24.05 & 3795 \\ \hline
\end{tabular}


\begin{figure*}[ht]
{\small \vspace{-90pt}
\begin{center}
		\begin{tabular}{@{\hspace{-22pt}}c@{\hspace{-17pt}}c @{\hspace{-17pt}}c@{\hspace{-17pt}}c} 
			\includegraphics[scale=0.32]{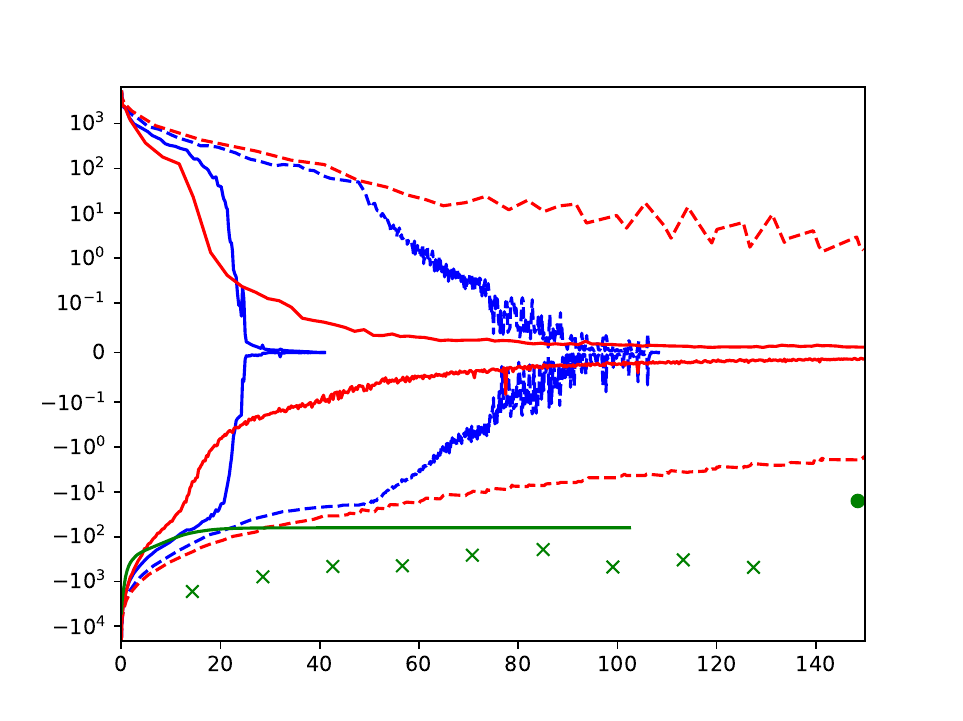} \begin{picture}(0, 0) \put(-85,90){\footnotesize 1CKK} \end{picture} &
			\includegraphics[scale=0.32]{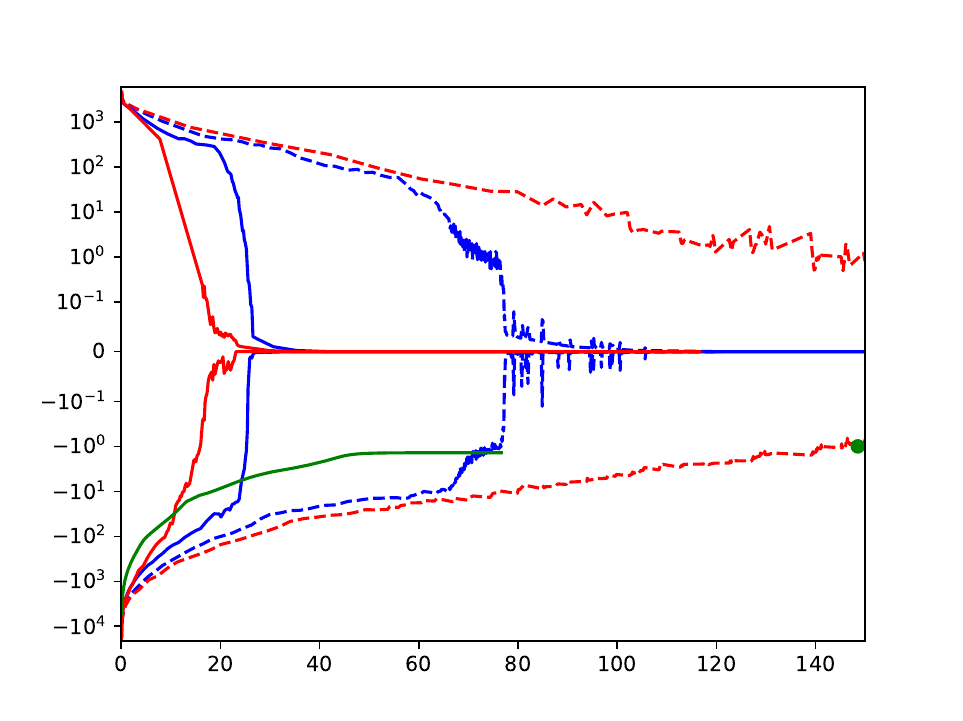} \begin{picture}(0, 0) \put(-85,90){\footnotesize 1CM1} \end{picture} &
			\includegraphics[scale=0.32]{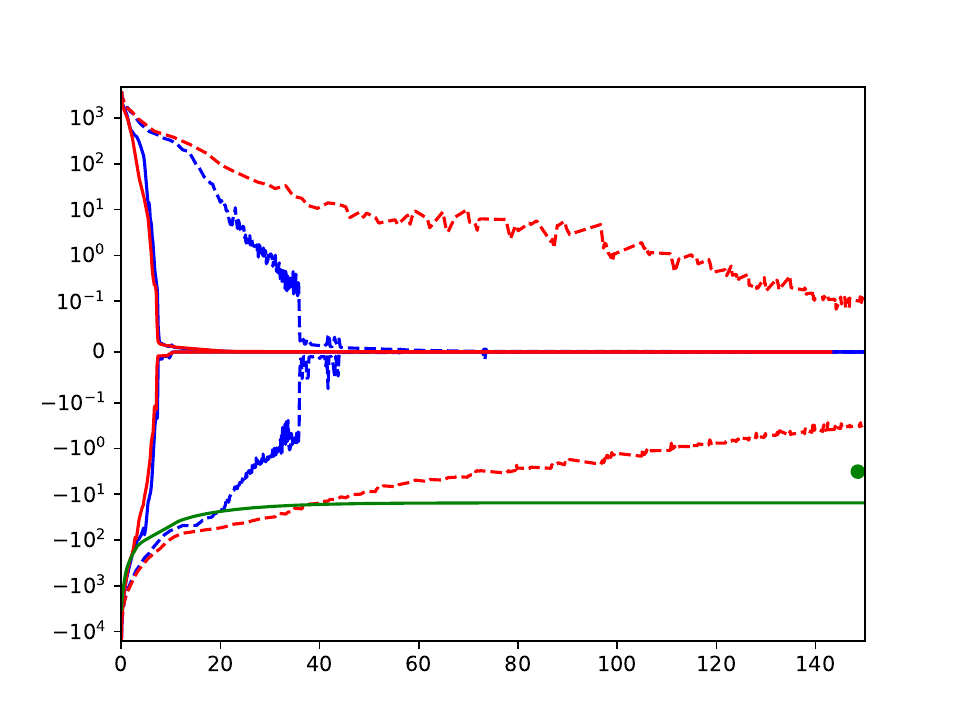} \begin{picture}(0, 0) \put(-85,90){\footnotesize 1SY9} \end{picture} &
			\includegraphics[scale=0.32]{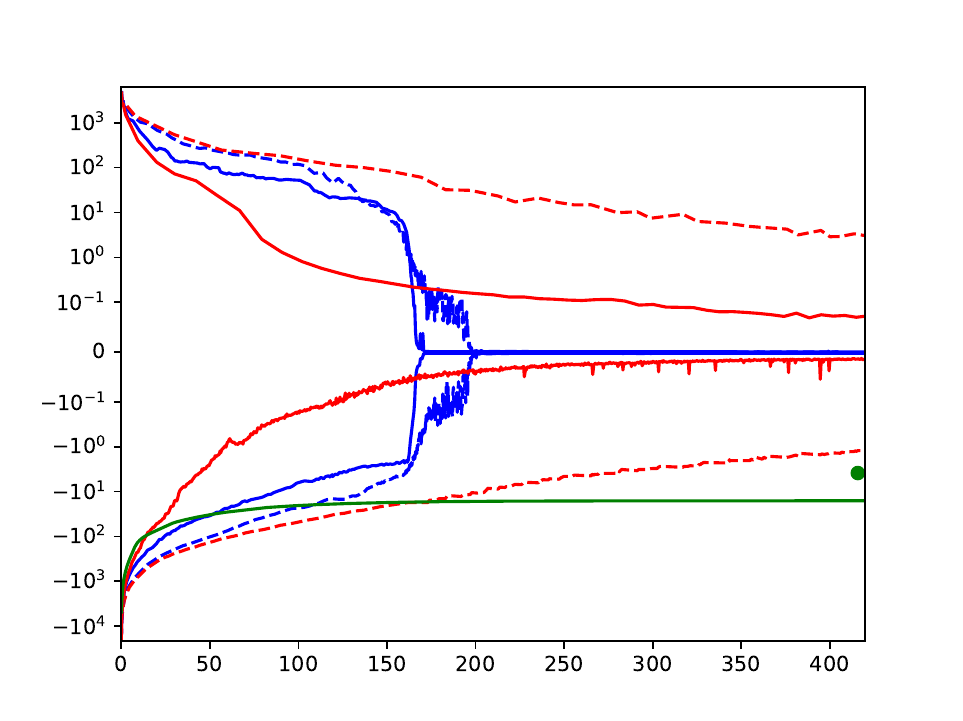} \begin{picture}(0, 0) \put(-85,90){\footnotesize 2BBN} \end{picture} 
		\end{tabular}\vspace{-5pt}

		\begin{tabular}{@{\hspace{-22pt}}c@{\hspace{-17pt}}c @{\hspace{-17pt}}c@{\hspace{-17pt}}c} 
			\includegraphics[scale=0.32]{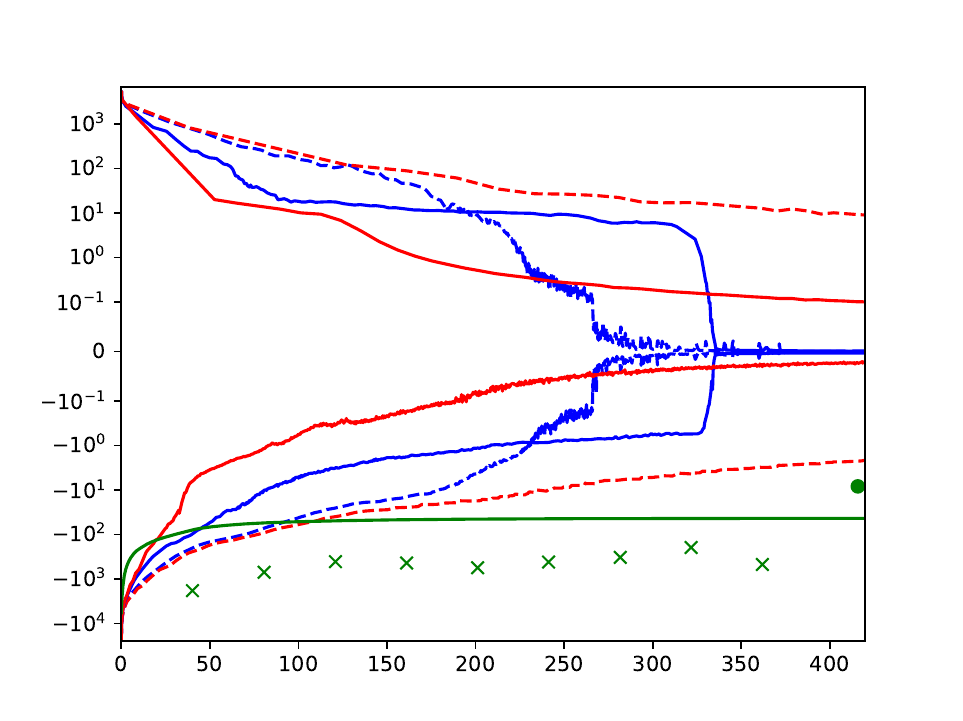} \begin{picture}(0, 0) \put(-85,90){\footnotesize 2BCX} \end{picture} &
			\includegraphics[scale=0.32]{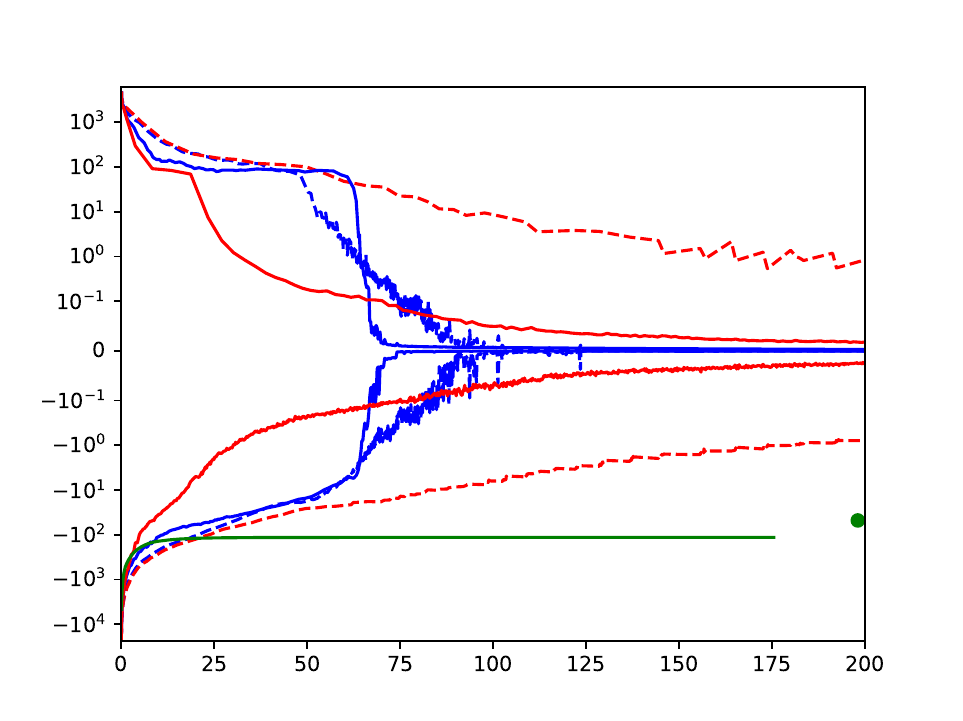} \begin{picture}(0, 0) \put(-85,90){\footnotesize 2BE6} \end{picture} &
			\includegraphics[scale=0.32]{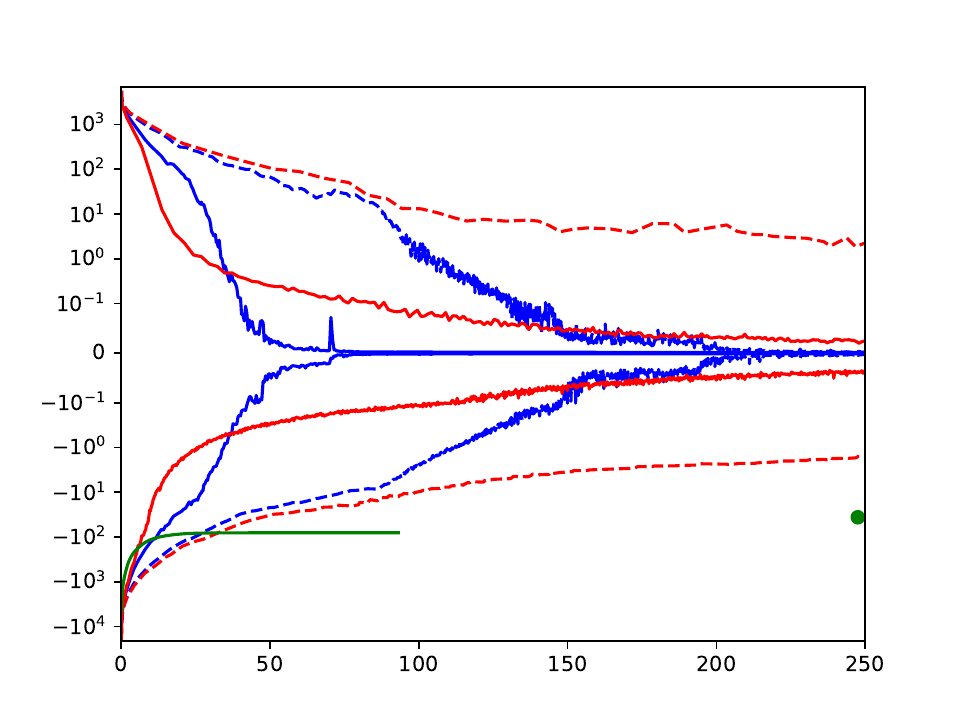} \begin{picture}(0, 0) \put(-85,90){\footnotesize 2F3Y} \end{picture} &
			\includegraphics[scale=0.32]{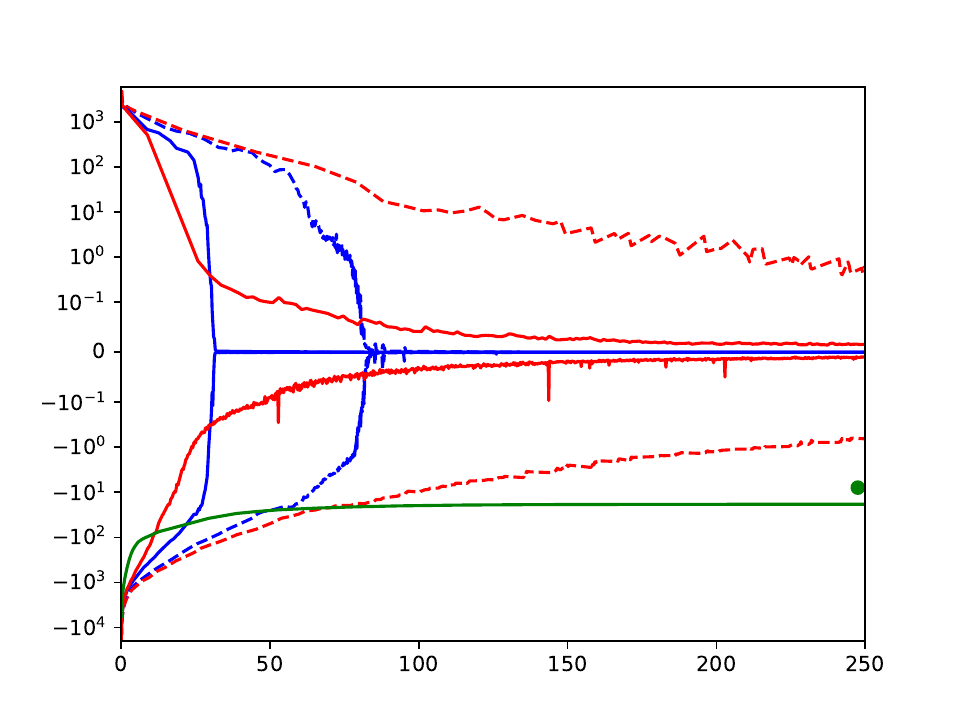} \begin{picture}(0, 0) \put(-85,90){\footnotesize 2FOT} \end{picture} 
		\end{tabular}\vspace{-5pt}
		
		\begin{tabular}{@{\hspace{-22pt}}c@{\hspace{-17pt}}c @{\hspace{-17pt}}c@{\hspace{-17pt}}c} 
			\includegraphics[scale=0.32]{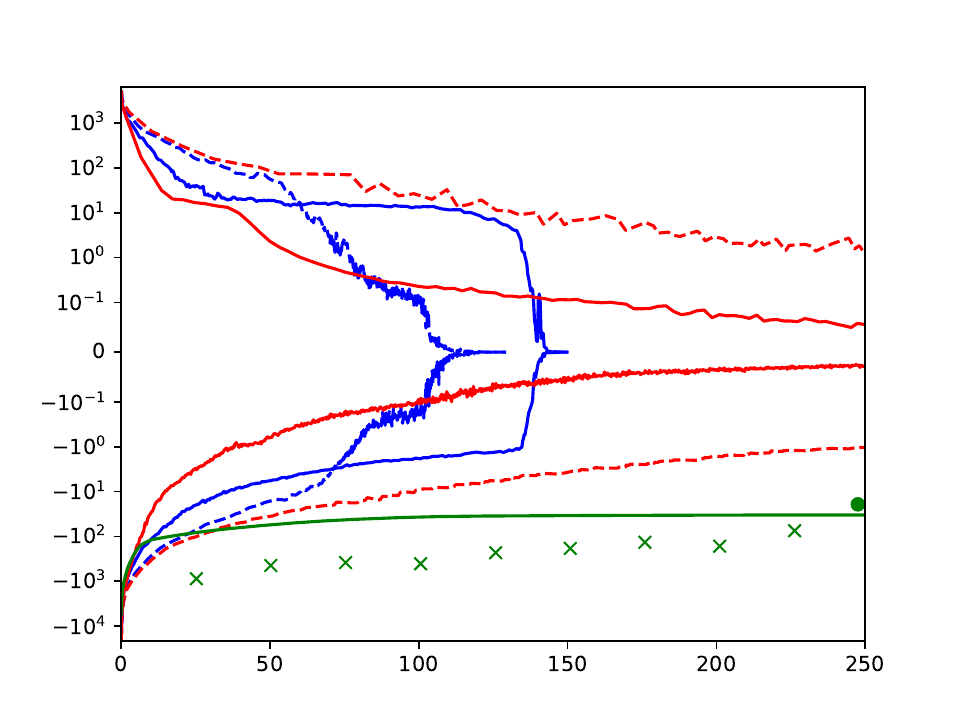} \begin{picture}(0, 0) \put(-85,90){\footnotesize 2HQW} \end{picture} &
			\includegraphics[scale=0.32]{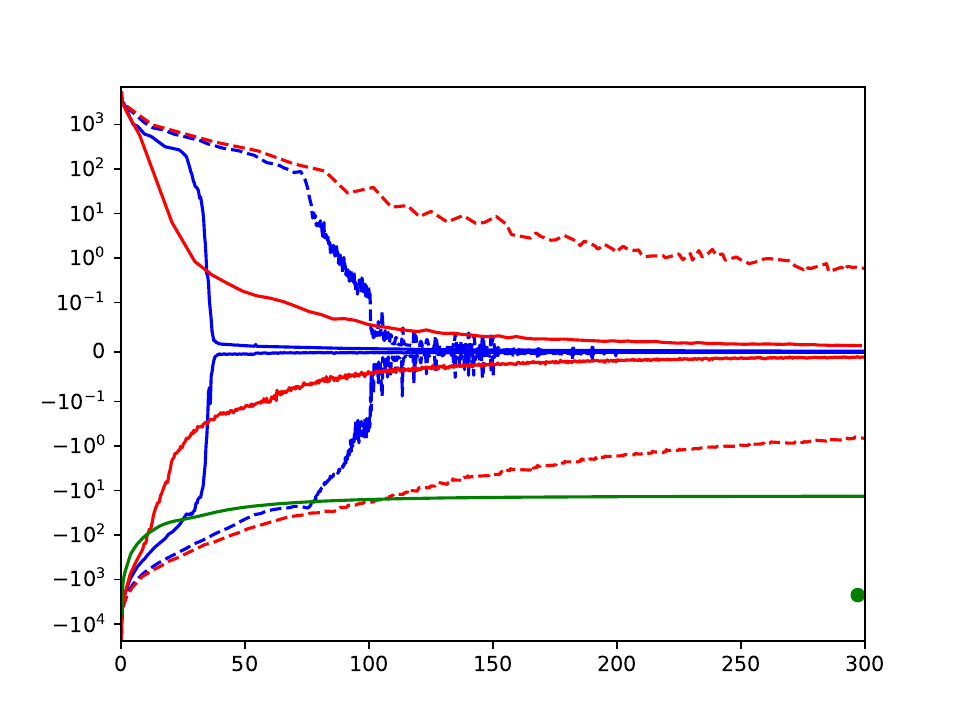} \begin{picture}(0, 0) \put(-85,90){\footnotesize 2O60} \end{picture} &
			\includegraphics[scale=0.32]{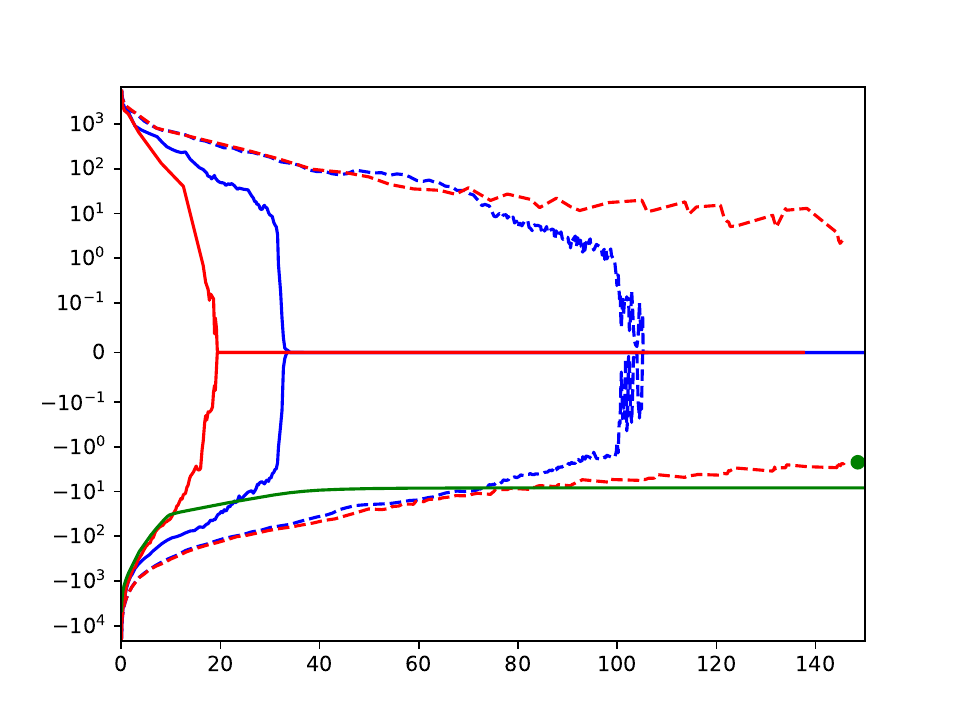} \begin{picture}(0, 0) \put(-85,90){\footnotesize 3BXL} \end{picture} &
			\includegraphics[scale=0.32]{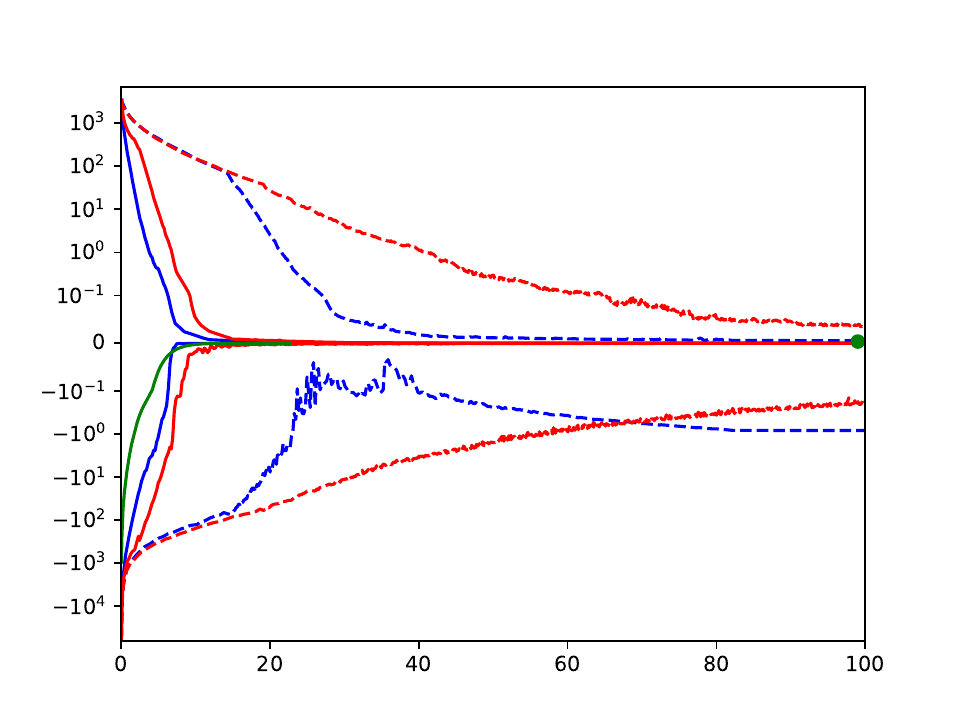} \begin{picture}(0, 0) \put(-85,90){\footnotesize pdb1b25} \end{picture} 
		\end{tabular}\vspace{-5pt}

		\begin{tabular}{@{\hspace{-22pt}}c@{\hspace{-17pt}}c @{\hspace{-17pt}}c@{\hspace{-17pt}}c} 
			\includegraphics[scale=0.32]{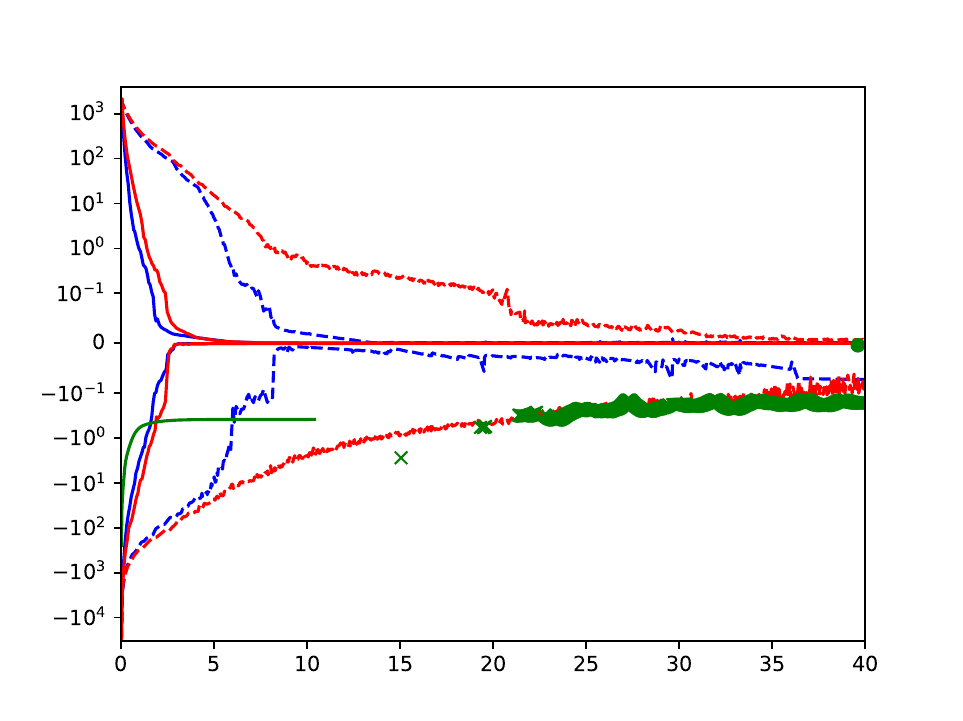} \begin{picture}(0, 0) \put(-85,90){\footnotesize pdb1d2e} \end{picture} &
			\includegraphics[scale=0.32]{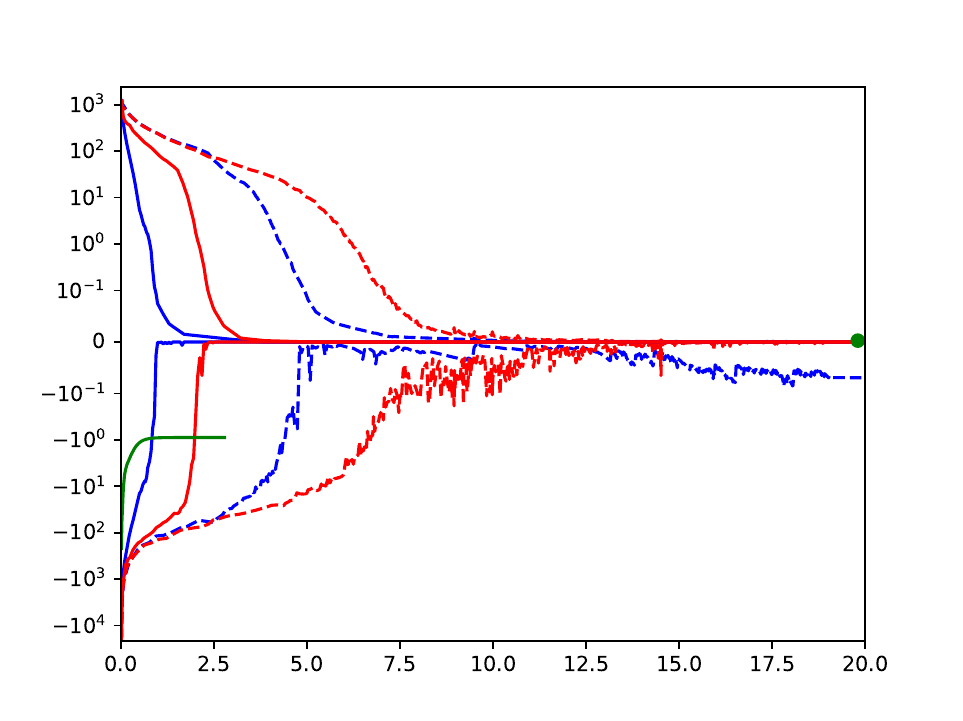} \begin{picture}(0, 0) \put(-85,90){\footnotesize pdb1fmj} \end{picture} &
			\includegraphics[scale=0.32]{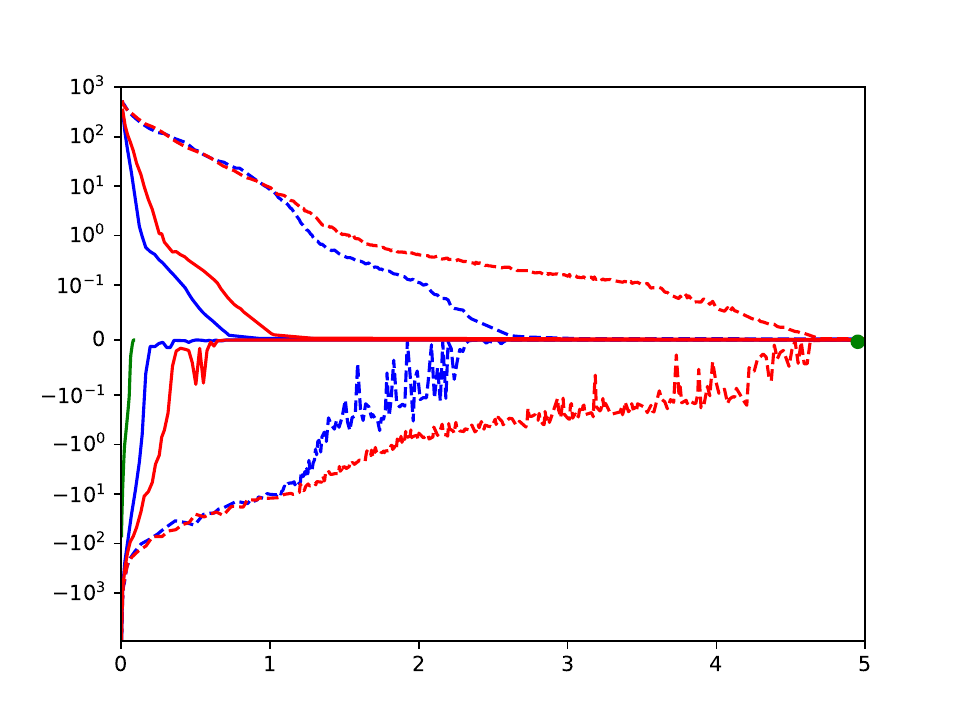} \begin{picture}(0, 0) \put(-85,90){\footnotesize pdb1i24} \end{picture} &
			\includegraphics[scale=0.32]{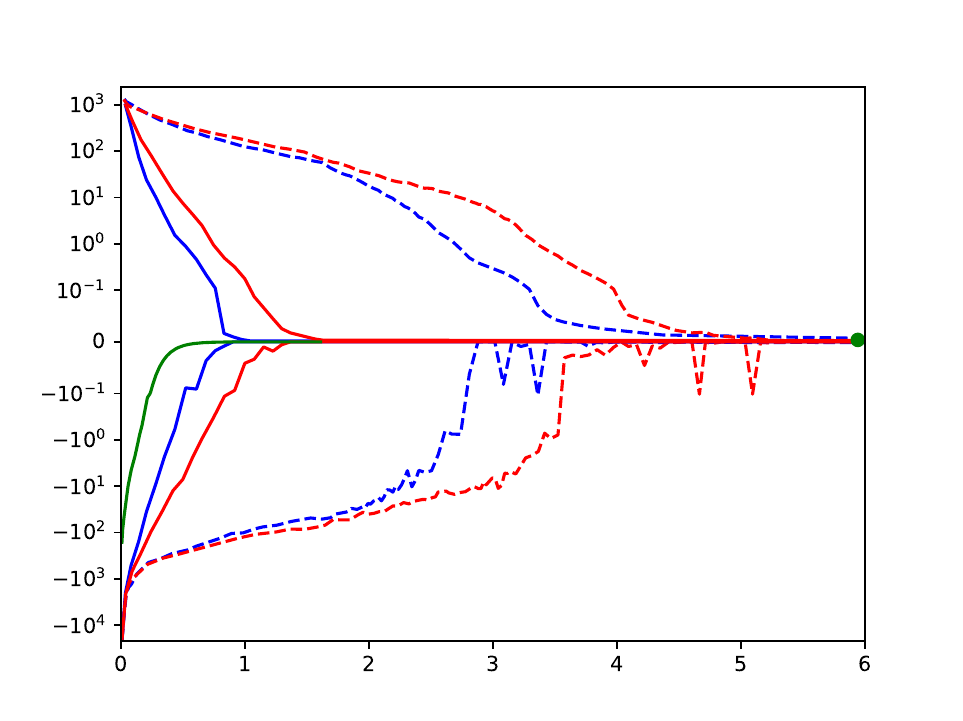} \begin{picture}(0, 0) \put(-85,90){\footnotesize pdb1iqc} \end{picture} 
		\end{tabular}\vspace{-5pt}\vspace{-7pt}


	\end{center}
	}
	\caption{Results on `protein-folding' ($\gamma=100$).\vspace{-7pt}
        }
\label{fig:protein-folding}
\end{figure*}

We tested 4 versions of FW algorithms that we term as ${\tt FW},{\tt FW}_{\tt conv},{\tt FW}^\ast,{\tt FW}^\ast_{\tt conv}$.
Subscript ${\tt conv}$ indicates option ${\tt conv}=\true$, while  ``$^*$'' means that in-face FW directions are used.
The FW method depends on parameter $\gamma$; values that are too small or too large result in slow convergence (see~\cite{Swoboda:CVPR19}).
For each family of problems we used the same $\gamma$ of the form $\gamma=10^k$, $k\in\mathbb Z$ with the best performance on one of the instances.
Since the ratios between allowed $\gamma$'s are rather large, we believe that it should be feasible to learn
such $\gamma$ (or even fine-tune it) for a given application using parameters such as the number of subproblems, etc.

In Fig.~\ref{fig:plots} and~\ref{fig:protein-folding} we plot lower and upper bounds as functions of time.
We added a constant to all values so that the average of the best known lower and upper bounds is zero,
and then used symmetric log scaling (`symlog' in python).
Note that the time for computing the upper bound (via OT) was not counted. The result of ADSAL is a single point
that we copied from~\cite{OpenGMBenchmark}. The runtime of ADSAL was always larger than the maximum $X$-range in the plots (sometimes significantly),
even though we used a slower machine: Intel Core i5-10210U CPU @ 1.60GHz and 16Gb RAM
vs. Intel Core i5-4570 CPU @ 3.20GHz and 32GB RAM used in~\cite{OpenGMBenchmark}.

\myparagraph{Discussion} We can see from the plots that ${\tt FW}^\ast_\sigma$ significantly outperforms ${\tt FW}_\sigma$ on `mrf-photomontage' and `protein-folding'
(both for $\sigma={\tt conv}$ and empty $\sigma$), as well as TRW-S and ADSAL.
\footnote{
As a single exception, ${\tt FW}$ starts giving better upper bound than ${\tt FW}^\ast$ once the lower bound reaches the optimal value (after roughly 500 seconds).
We do not have an explanation for this behavior. Perhaps, an alternative method for extracting a primal solution should be used in this regime.} 
This suggests that FW with in-face FW directions is the current state-of-the-art LP solver for these applications.
On two other applications, however (`matching' and `object-seg') ${\tt FW}^\ast_\sigma$ and ${\tt FW}_\sigma$ are roughly similar
and in general outperformed by other techniques, e.g.\ AD3 \& Gurobi on `matching' and TRW-S on `object-seg'.

On most plots upper and lower bounds seem to be converging to each other. We conclude that the FW approach
can be used for finding both primal and dual solutions of relaxation~\eqref{eq:BLP}.


{\small
\bibliographystyle{ieee_fullname}
\bibliography{FW}
}

\end{document}